\documentclass[final]{article}

\newcommand{\modif}[1]{#1}

\usepackage{amsmath, amstext, amssymb}
\usepackage{units}
\usepackage{subfigure}
\usepackage{float}
\usepackage{multirow}
\def\paragraph#1{{\bf #1\ }}

\def\NN{\mbox{\rm I\hspace{-0.50ex}N} }
 \def\OO{\rm \hbox{O\kern-.34em\raise.47ex
         \hbox{$\scriptscriptstyle |$}\kern-.46em\raise.47ex
         \hbox{$\scriptscriptstyle |$}\kern+0.5 em }}
\def\RR{\mbox{\mathrm I\hspace{-0.50ex}R} }

\def\hcboxcm#1#2{\hbox to #1{\hfill #2 \hfill}}

\def\null{\hbox{}}
\let\eps\varepsilon

\def\tn1{\widetilde n_1}
\def\tn2{\widetilde n_2}
\def\tn{\widetilde n }

\def\be{\begin{equation}}
\def\ee{\end{equation}}
\def\bea{\begin{eqnarray}}
\def\eea{\end{eqnarray}}
\def\bean{\begin{eqnarray*}}
\def\eean{\end{eqnarray*}}

\def\RR{{\mathrm{ I~\hspace{-1.15ex}R}}}
\def\NN{{\mathrm{ I~\hspace{-1.15ex}N}}}
\def\qquad{{\quad\quad}}

\def\={\, = \, }

\def\Box{\leavevmode\vbox{\hrule
     \hbox{\vrule\kern4pt\vbox{\kern4pt}%
           \vrule}\hrule}}
\def\blackbox{\leavevmode\vrule height 5pt width 4pt depth 0pt\relax}
\catcode`@=11

\def\eqalign#1{\null\,\vcenter{\openup1\jot \m@th
   \ialign{\strut \hfil$\displaystyle{##}$ & $\displaystyle{{}##}$\hfil
      \crcr#1\crcr}}\,}
\def\eqalignrll#1{\null\,\vcenter{\openup1\jot \m@th
   \ialign{\strut \hfil$\displaystyle{##}$ & $\displaystyle{{}##}$\hfil
    & $\displaystyle{{}##}$\hfil
      \crcr#1\crcr}}\,}
\def\eqalignrcl#1{\null\,\vcenter{\openup1\jot \m@th
   \ialign{\strut \hfil$\displaystyle{##}$ &\hfil $\displaystyle{{}##}$\hfil
    & $\displaystyle{{}##}$\hfil
      \crcr#1\crcr}}\,}
\def\eqalignno#1{\displ@y \tabskip\@centering
  \halign to\displaywidth{\hfil$\@lign\displaystyle{##}$\tabskip\z@skip
    &$\@lign\displaystyle{{}##}$\hfil\tabskip\@centering
    &\llap{$\@lign##$}\tabskip\z@skip\crcr
    #1\crcr}}
\newcounter{appendix}
\newcounter{sectionz}
\def\appendix{\advance\c@appendix by 1
   \def\thesectionz {\Alph{appendix}}
\def\thesection{\Alph{appendix}} 
   \ifnum\c@appendix=1 \setcounter{section}{-1} \fi
   \@startsection {section}{1}{\z@}{-3.5ex plus -1ex minus 
  -.2ex}{2.3ex plus .2ex}{\large\bf}}

\catcode`@=12
 \newtheorem{lemme}{Lemma}[section]  
 \newtheorem{theorem}[lemme]{Theorem}
 \newtheorem{remark}[lemme]{Remark} 

\def\deblem{\begin{lemme} }
\def\finlem{\end{lemme}}
\def\debthm{\begin{theorem} }
\def\finthm{\end{theorem}}
\def\debprop{\begin{proposition}}
\def\finprop{\end{proposition}}
\def\debcor{\begin{corollary}}
\def\fincor{\end{corollary}}
\def\debdef{\begin{definition}}
\def\findef{\end{definition}}
\def\debrem{\begin{remark}}
\def\finrem{\null\hfill\blackbox\end{remark}}
\def\debproof{{\bf Proof: \ }}
\def\finproof{\null\hfill {$\blackbox$}\bigskip}  
\def\NN{\mathbb{N}}
\def\OO{\mathbb{O}}

\def\RR{\mathbb{R}}

\usepackage{amsfonts}
\usepackage{amssymb}
\usepackage{amsbsy}

\usepackage{stmaryrd}
\usepackage[dvips]{graphicx}
\usepackage{psfrag}
\usepackage[hang,center]{caption}
\usepackage{float}
\usepackage{cite}

\usepackage{color}
\usepackage{cases}

\title{An asymptotic preserving scheme for strongly anisotropic elliptic problems}

\author{Pierre Degond\footnotemark[2]\ \footnotemark[3] \and  Fabrice Deluzet\footnotemark[2]\ \footnotemark[3] \and Claudia Negulescu\footnotemark[4] }

\begin{document}
\maketitle

\renewcommand{\thefootnote}{\fnsymbol{footnote}}

\footnotetext[2]{Universit\'e de Toulouse, UPS, INSA, UT1, UTM, Institut de Math\'ematiques de Toulouse, F-31062 Toulouse, France}
\footnotetext[3]{CNRS, Institut de Math\'ematiques de Toulouse UMR 5219, F-31062 Toulouse, France}
\footnotetext[4]{CMI/LATP, Universit\'e de Provence, 39 rue Fr\'ed\'eric Joliot-Curie 13453 Marseille cedex 13} 

\renewcommand{\thefootnote}{\arabic{footnote}}

\begin{abstract}
In this article we introduce an asymptotic preserving scheme designed
to compute the solution of a two dimensional elliptic equation
presenting large anisotropies. We focus on an anisotropy aligned with one direction, the dominant part of the elliptic operator being supplemented with Neumann boundary conditions. A new scheme is introduced which allows an accurate resolution of this elliptic equation for an arbitrary anisotropy ratio.
\end{abstract}



\pagestyle{myheadings}
\thispagestyle{plain}
\markboth{P. DEGOND, F. DELUZET  AND C. NEGULESCU}{An asymptotic preserving scheme for strongly anisotropic elliptic problems}

\section{Introduction}
The objective of this paper is to introduce an efficient and accurate
numerical scheme to solve a strongly anisotropic elliptic problem of the
form
\be \label{Aniso}
\left\{
\begin{array}{l}
\displaystyle  - \nabla \cdot \left ( \mathbb{A} \nabla \phi \right ) = f \,, \quad \textrm{in} \,\, \Omega\\[3mm]
\displaystyle \phi = 0 \quad \textrm{on} \,\, \partial \Omega_D\,, \quad \displaystyle \partial_z \phi = 0 \quad  \textrm{on} \,\, \partial \Omega_z\,,
\end{array}
\right.
\ee
where $\Omega \subset \RR^2$ or $\Omega \subset \RR^3$ is a domain, with boundary $\partial \Omega=\partial \Omega_D \cup \partial\Omega_z$ and the diffusion matrix $\mathbb{A}$ is given by
$$
\mathbb{A} =\left(
  \begin{array}[c]{cc}
    A_\perp & 0 \\
    0 & \frac{1}{\varepsilon}A_z 
  \end{array}\right) \,.
$$
The terms $A_\perp$ and $A_z$ are of the same order of magnitude,
whereas the parameter $0<\eps<1$ can be very small, provoking thus the
high anisotropy of the problem. In the present paper the considered anisotropy direction is fixed and is aligned with the $z$-axis of a Cartesian coordinate system. The method presented here is extended in some forthcoming works to more general anisotropies \cite{brull}.

\noindent Anisotropic problems are common in mathematical modeling and
numerical simulation. Indeed they occur in several fields of
applications such as flows in porous media~\cite{porous1,TomHou},
semiconductor modelling \cite{semicond}, quasi-neutral plasma
simulations \cite{Navoret}, image processing \cite{Weickert,image1},
atmospheric or oceanic flows \cite{ocean}, and so on, the list being
not exhaustive.  More specifically high anisotropy aligned with one direction may occur in shell problems or simulation in stretched media. The initial motivation for the present work is closely related to magnetized plasma simulations such as atmospheric \cite{Kelley2, Kes_Oss} or inertial fusion plasmas \cite{Beer, Sangam} or plasma thrusters \cite{SPT}. In this context, the medium is structured by the magnetic field. Indeed, the motion of charged particles in planes perpendicular to the magnetic field is governed by a fast gyration around the magnetic field lines. This explains the large number of collisions the particles encounter in the perpendicular plane, whereas the dynamic in the parallel direction is rather undisturbed. As a consequence the particle mobilities in the perpendicular and parallel directions differ by many orders of magnitude. In the context of ionospheric plasma modelling \cite{BCDDGT_4,Hysell}, the ratio of the aligned and transverse mobilities (denoted in this paper by $\eps^{-1}$)  can be as huge as ten to the power ten. The relevant boundary conditions in many fields of application are periodic (for instance in simulations of tokamak plasmas on a torus) or Neumann boundary conditions (see for instance \cite{BCDDGT_1} for atmospheric plasmas). The system \eqref{Aniso} is thus a good model to elaborate a robust numerical method. \\

\noindent The main difficulties with the resolution of problem (\ref{Aniso}) are of numerical nature, as solving this singular perturbation problem for small $0<\eps \ll1$ is rather delicate. Indeed, replacing in the anisotropic elliptic equation $\eps$ by zero, yields an ill-posed problem, which has an infinite number of solutions (namely all functions which are constant in the $z$-direction). This feature is translated in the discrete case (after the discretization of the problem) into  a linear system which is very ill-conditioned for $\eps \ll1$, due to the different order of
magnitudes of the various terms.  As a consequence standard numerical methods for the resolution of linear systems lead to important numerical costs and unacceptable numerical errors.

  \modif{More generally, this numerical difficulty arises when the boundary
  conditions supplied to the dominant $O(1/\varepsilon)$ operator lead to an ill-posed problem with a multiplicity of solutions. This is the case for Neumann boundary conditions, but also of periodic boundary conditions. If instead, the boundary conditions are such that the dominant operator gives a well-posed problem with a unique solution, this difficulty
  vanishes as the leading operator alone will suffice to completely determine the limit solution. In this case, one can resort to standard methods. This is the case of Dirichlet or Robin boundary conditions. In spite of the fact that the problem addressed in the present paper arises only with specific boundary conditions, it has a considerable impact in many physics problem, such as  plasmas, geophysical flows, plate and shells, etc. In this paper, we will focus on Neumann boundary conditions because they represent a larger range of  physical applications, but we could address periodic boundary conditions in a similar way. }

\modif{Numerical methods for anisotropic elliptic problems have been
  extensively investigated in the literature.}
 Depending on the underlying
physics, distinct numerical methods are developed. For example domain
decomposition (Schur complement) and multigrid techniques, using
multiple coarse grid corrections are adapted to anisotropic equations
in \cite{Giraud,Koronskij} and \cite{Gee,Notay}. For anisotropy
aligned with one (or two directions), point (or plane) smoothers are shown to be very efficient \cite{ICASE}. A problem very similar to \eqref{Aniso} is addressed in \cite{Vladimir}, treated via a parametrisation technique, and seems to give good results for rather large anisotropy ratios. However, these techniques are only developed in the context of an elliptic operator with a dominant part supplemented with Dirichlet boundary conditions.\\
\modif{An alternative approach} 
for dealing with highly anisotropic problems is based on a mathematical reformulation of the continuous problem, in order to obtain a more harmless problem, which can be solved numerically in an uncomplicated manner. In this category can be situated for example asymptotic models, describing for small values of the asymptotic parameter $\varepsilon$ the evolution of an approximation $\tilde{\phi}$ of the solution of (\ref{Aniso}) \cite{BCDDGT_1,Keskinen}. However, these asymptotic models are precise only for
$\eps \ll 1$, and cannot be used on the whole range of values covered by the physical parameter $\varepsilon$. Thus model coupling methods have to be employed. In sub-domains where the limit model is no longer valid, the original model has to be used, which means that a model coupling strategy has to be developed. However the coupling strategy requires the existence of an area where both models are valid and still demands an accurate numerical method for the resolution of the original model (i.e. the anisotropic elliptic problem) with large anisotropies. This can be rather undesirable. \\

\noindent In this paper, we present an original numerical algorithm belonging to the second approach. A reformulation of the continuous problem (\ref{Aniso}) will permit us to solve this problem in an inexpensive way and accurately
enough, independently of the parameter $\varepsilon$. This scheme is related to the {\it Asymptotic Preserving} numerical method introduced in \cite{ShiJin}. These techniques are designed to provide computations in various regimes without any restriction on the discretization meshes and with the additional property to converge towards the solution of the limit problem when the asymptotic parameter goes to zero. The derivation of such Asymptotic Preserving methods requires first the identification of the limit model. 
For singular perturbation problems, a reformulation of the problem is required in order to derive a set of equations containing both the initial and the limit model with a continuous transition from one regime to another, according to the values of the parameter $\varepsilon$. This reformulated system of equations sets the foundation of the AP-scheme. Other singular perturbations have already been explored in previous studies, for instance quasi-neutral or gyro-fluid limits \cite{Crispel, Sangam}. These techniques have been first introduced for non-stationary systems of equations, for which the time discretization must be studied with care in order to guarantee the asymptotic preserving property. For the anisotropic elliptic equation investigated in this article, we only need to precise the reformulated system and provide a discretization of this one.\\

\noindent The outline of this paper is the following. Section~\ref{SEC2} of this article presents first the initial anisotropic elliptic model. In the remainder of this paper, it will be referred to as the Singular-{\it{Perturbation}} model (P-model). The reformulated system (referred to as the {\it Asymptotic Preserving} formulation or  AP-formulation) is then derived. It relates on a decomposition of the solution $\phi(x,z)$ according to its {\it mean part} $\bar \phi(x)$ along the $z$ coordinate and a {\it fluctuation} $\phi'(x,z)$ consisting of a correction to the mean part needed to recover the full solution.  The mean part $\bar {\phi}(x)$ is solution of an $\eps$-independent
elliptic problem, and the fluctuation $\phi'(x,z)=\phi(x,z)-\bar{\phi}(x)$ is given by a well-posed $\eps$-dependent elliptic problem. The
advantage is that the $\eps$-dependent problem for the fluctuation
is well-posed and solvable in an inexpensive way, and this uniformly
in $\varepsilon$. In the limit $\eps \rightarrow 0$ the AP-formulation
reduces to the so called {\it Limit} model ($L$-model), whose solution
is an acceptable approximation of the P-model solution  for
$\eps \ll1$. The present derivation is carried out in the framework of
an anisotropy aligned along one  axis of a Cartesian coordinate
system. In the context of magnetized plasma simulations, this initial
work is extended in a forthcoming work for the three dimensional
case in curvilinear coordinates, designed to fit a more complex
magnetic field topology (i.e. anisotropy direction)
\cite{BCDDGT_4}. The main constraints of this method reside in the
construction of the mean part which necessitates the integration of
the solution along the anisotropy direction. This operation is easily carried out in the context of coordinates adapted with the anisotropy direction. However,  an extension of the techniques presented here is currently developed for non-adapted coordinates
\cite{brull}. \\
\noindent Section \ref{SEC3} is devoted to 
 the numerical implementation of the AP-formulation. Numerical results are
 then presented for a test case, and the three approaches (AP-formulation, straight discretization and resolution of the P-model and L-model) are compared
 according to the precision of the approximation for different values of $\eps$. In section \ref{SEC4} we shall rigorously analyse the
 convergence of the AP-scheme. Error estimates will be established
 which underline the advantages of the AP-scheme as compared to the
 initial Singular Perturbation model and the Limit model.
 
 \modif{Current research directions are concerned with the adaptation of the present technique to the case of arbitrary spatially varying anisotropies, without adaptation of the coordinate system to the direction of the anisotropy. These developments will allow the treatment on nonlinear problems, when the diffusion tensor (and its principal directions) depend on the solution itself. This treatment will involve iterative methods which, at each iterate, will reduce the problem to the solution of a linear anisotropic diffusion problem. }
\section{The asymptotic preserving formulation} \label{SEC2}
For simplicity we shall consider in this paper the two-dimensional
problem, posed on a rectangular domain $\Omega= \Omega_x \times
\Omega_z$, where $\Omega_x \subset \RR$ and $\Omega_z \subset \RR$ are
intervals. The ideas exposed here can be extended without any problems
to the more physical three-dimensional domain, with two transverse
directions $(x,y)$ and an anisotropy direction aligned with the $z$-direction. 
In this section we introduce the Singular Perturbation Model, the
Limit Model and the Asymptotic Preserving formulation.
\subsection{The Singular Perturbation Model (P-model)}\label{SEC21}
The main concern of this paper is the numerical resolution of the
following anisotropic, elliptic problem, called in the sequel Singular Perturbation Model
\begin{equation}\label{eq:ellipti:original}
(P)\,\,\,
\left\{ 
\begin{array}{l}
 \displaystyle - \nabla \cdot \left ( \mathbb{A} \nabla \phi \right ) = f\,,\quad \text{in} \quad \Omega\,, \\[3mm]
 \displaystyle  \frac{\partial \phi}{\partial z} = 0 \quad \text{on} \quad
\Omega_x \times \partial \Omega_z\,,\qquad
 \displaystyle \phi = 0 \quad \text{on} \quad \partial \Omega_x \times \Omega_z \,.
\end{array}
\right.
\end{equation}
The anisotropy of the media is modeled  via the definition of the diffusion matrix $\mathbb{A}$
\begin{equation}
  \mathbb{A} =\left(
  \begin{array}[c]{cc}
    A_\perp & 0 \\
    0 & \frac{1}{\varepsilon}A_z 
  \end{array}\right) \,,
\end{equation}
where $A_\perp(x,z)$ and $A_z(x,z)$ are given functions with
comparable order of magnitudes. The source term $f(x,z)$ is given and the parameter $\varepsilon$ is small
compared to both $A_\perp$ as well as $A_z$. The medium becomes more
anisotropic as the value of  $\varepsilon$ goes to zero. 

\subsection{The limit regime (L-model)}\label{SEC22}
In this section we establish that in the limit $\eps \rightarrow 0$ the solution of the perturbation model converges towards $\bar\phi$, solution of the L-model defined by 
\begin{equation}\label{eq:ellipti:S}
(L)\,\,\,
\left\{
\begin{array}{l}
\displaystyle- \frac{\partial}{\partial x} \left( \bar{A}_\perp \frac{\partial \bar{\phi}}{\partial x}\right)   = \bar{f}(x) \,,\quad \text{in} \quad \Omega_x\,, \\[3mm]
 \displaystyle \bar{\phi} = 0 \quad \text{on} \quad \partial \Omega_x \,,
\end{array}
\right.
\end{equation}
where overlined quantities designate averages over the z-coordinate~:
$$\bar f(x) = \modif{\frac{1}{|\Omega_z|}\int_{\Omega_z}} f(x,z) \, dz.$$

\noindent First we can rewrite the P-model as 
\begin{equation}\label{eq:ellipti:original_bis}
(P)\,\,\,
\left\{
\begin{array}{l}
\displaystyle- \frac{\partial}{\partial x} \left( {A}_\perp \frac{\partial {\phi}}{\partial x}\right)  - \frac{1}{\eps} \frac{\partial}{\partial z} \left( {A}_z \frac{\partial {\phi}}{\partial z}\right) = {f} \,,\quad \text{in} \quad \Omega\,, \\[3mm]
 \displaystyle \frac{\partial \phi}{\partial z} = 0 \quad \text{on} \quad \Omega_x  \times \partial \Omega_z \,,\quad \phi = 0 \quad \text{on} \quad \partial \Omega_x  \times \Omega_z \,,
\end{array}
\right.
\end{equation}
and integrating along the $z$-coordinate gives
  \begin{equation}\label{aniso:integrated}
    \frac{\partial}{\partial x} \left( \overline{A_\perp \frac{\partial \phi}{\partial x} } \right)  = \bar f(x) \,.
  \end{equation}
 This equation holds for any $\eps > 0$. Now, letting formally
 $\eps$ tend to zero in \eqref{eq:ellipti:original_bis} yields
 the reduced model (R-model) 
\begin{equation}\label{redu}
(R)\,\,\,
\left\{
\begin{array}{l}
\displaystyle - \frac{\partial }{\partial z} \left(A_z \frac{\partial \phi}{\partial z} \right)  = 0 \,,\quad \text{in} \quad  \Omega\,, \\[3mm]
\displaystyle \frac{\partial \phi}{\partial z} = 0 \quad \text{on}\quad
\Omega_x \times \partial \Omega_z \,, \qquad
\displaystyle \phi = 0 \quad \text{on} \quad \partial \Omega_x \times \Omega_z\,.
\end{array}
\right.
\end{equation}
The functions verifying this ill-posed R-model are constant
  along the $z$-coordinate. Thus including this asymptotic limit
  property into equation~\eqref{aniso:integrated} gives rise to the L-model~\eqref{eq:ellipti:S}, verified by the solution of the Singular Perturbation model in the limit $\eps \rightarrow 0$.

\begin{remark}
The L-model is the singular limit of the original P-model~\eqref{eq:ellipti:original}. It provides an accurate approximation of the P-solution only for small values of $\eps$. The P-model is valid for all $0<\eps <1$, but numerically impracticable for $\eps\ll 1$. Indeed working with a finite precision, the asymptotic model degenerates into the R-model defined by \eqref{redu} as $\eps$ vanishes. This R-model is ill-posed since it exhibits an infinite amount of solutions $\phi = \tilde \phi(x)$, depending only on the variable $x$. This implies that the discretization matrix derived from the P-model is very
ill-conditioned for small $0<\eps \ll 1$. This point is addressed by the numerical experiments of section~\ref{SEC33}.
 Consequently, in a domain where $\eps$ varies significantly, a model
 coupling method has to be developed in order to exploit the validity
 of each model, the P- and L-model. This can be rather undesirable. In the next section we
shall present an alternative approach, which is based on a reformulation
of the Singular-Perturbation model providing a means of computing an accurate numerical approximation of the solution for all values $0<\eps<1$. 
\end{remark}
\begin{remark}
\modif{The asymptotics is totally different in the case of Dirichlet boundary conditions. In this case, the R-model is well posed, with a unique solution, and there is no difficulty anymore. Any standard numerical solution of the P-model will converge to that of the R-model. In other words, with Dirichlet boundary conditions, the perturbation becomes regular and the limit solution is fully determined by the formal limit system. The situation and the difficulty addressed in the present paper require that the R-model be ill-posed. This is the case with Neumann boundary conditions (which is the framework chosen here) but also with periodic boundary conditions, or any other boundary condition which would result in an ill-posed R-model. }  
\end{remark}

\subsection{The Asymptotic Preserving reformulation (AP-formulation)}\label{SEC23}
In order to circumvent the just described numerical difficulties in handling the Singular Perturbation model, we introduce a reformulation, which permits a transition from the initial $P$-model to its singular limit (L-model), as $\eps \rightarrow 0$.\\
For this, we shall decompose each quantity $f(x,z)$ into its mean value $\bar f(x)$
along the $z$ coordinate and a fluctuation part $f'(x,z)$. For
simplicity reasons let in the following $\Omega_x:=(0,L_x)$ and
$\Omega_z:=(0,L_z)$. Then 
\begin{align}\label{eq:def:decomp}
  f(x,z) &= \bar f(x) + f'(x,z) \,,
\end{align}
with 
\begin{align}\label{eq:def:decomp:bis}
  \bar f(x) := \frac{1}{L_z} \int_{0}^{L_z} f(x,z) dz \,, \quad f'(x,z) := f(x,z) - \bar f(x) \,.
\end{align}
Note that we have the following properties
\begin{align}
 & \bar{f'} = 0 \,,& \overline{ \left( \nicefrac{\partial f}{\partial x} \right) }= \nicefrac{\partial \overline{f}}{\partial x}\,,\qquad    & \displaystyle \overline{fg} = \bar f \bar g + \overline{f'g'} \,,\label{eq:prop:mean}\\
 & \nicefrac{\partial f}{\partial z} = \nicefrac{\partial f'}{\partial z} \,, & \left( \nicefrac{\partial f}{\partial x} \right)' = \partial \nicefrac{f'}{\partial x}\,,\qquad  & (fg)' = f'g' - \overline{f'g'}+ \bar f g' + f' \bar g\,. \label{eq:prop:fluct}
\end{align}
Taking now the mean of the elliptic equation
\eqref{eq:ellipti:original_bis} along the $z$-coordinate, we get thanks
to~\eqref{eq:prop:mean} and~\eqref{eq:prop:fluct}, an equation for the
evolution of the mean part $\bar{\phi}(x)$
\begin{equation}\label{eq:compute:mean}
(AP1)\,\,\,
\left\{
\begin{array}{l}
\displaystyle  - \frac{\partial}{\partial x}\left( \bar A_\perp
\frac{\partial \bar \phi}{\partial x} \right) = \bar f +
\frac{\partial}{\partial x}\left( \overline{ A_\perp' \frac{\partial
\phi'}{\partial x}} \right)\,,\quad \text{in} \quad  \Omega_x\,, \\[3mm]
\displaystyle  \overline{\phi} = 0 \quad \text{on} \quad \partial \Omega_x\,.
\end{array}
\right.
\end{equation}
Substracting from \eqref{eq:ellipti:original_bis} this mean
equation~\eqref{eq:compute:mean}, gives rise to the evolution equation
for the fluctuation part $\phi'(x,z)$
\begin{equation}\label{sys:phi:prime}
\hspace*{-0.2cm}(AP2)\,\,\,
\left\{
\begin{array}{l}
\displaystyle - \frac{\partial }{\partial z} \left(A_z \frac{\partial \phi'}{\partial z} \right) - \varepsilon \frac{\partial}{\partial x}\left(A_\perp \frac{\partial \phi'}{\partial x}  \right) + \varepsilon \frac{\partial}{\partial x}\left(\overline{A_{\perp}' \frac{\partial \phi'}{\partial x} }\right) = \\[3mm]
 \displaystyle \hspace{5.2cm} \varepsilon f' + \varepsilon \frac{\partial}{\partial x} \left(A_\perp' \frac{\partial \overline{\phi}}{\partial x} \right)\,,\quad \text{in} \,\,\,  \Omega\,, \\[3mm]
 \displaystyle \frac{\partial \phi'}{\partial z} = 0 \quad \text{on}\quad
\Omega_x \times \partial \Omega_z \,, \qquad
\displaystyle \phi' = 0 \quad \text{on} \quad \partial \Omega_x \times \Omega_z\,, \\[3mm]
 \displaystyle  \overline{ \phi'}= 0\,, \quad \text{in} \quad \Omega_x\,.
\end{array}
\right.
\end{equation}
Thus we have replaced the resolution of the initial Singular Perturbation model
(\ref{eq:ellipti:original_bis}) by the resolution of the system
(\ref{eq:compute:mean})-(\ref{sys:phi:prime}), which will be done
iteratively. Starting from a guess function $\phi'$, equation (\ref{eq:compute:mean}) gives the mean
value $\overline{\phi}(x)$, which inserted in (\ref{sys:phi:prime})
shall give the fluctuation part $\phi'(x,z)$ and so on.\\

The constraint $\overline{\phi'} =0$ in (\ref{sys:phi:prime})
(which is automatic for $\eps>0$, as explained in Remark
  \ref{rem1}) has the essential consequence that the conditioning of the discretized system becomes $\eps$-independent, because the problem (\ref{sys:phi:prime}) reduces in the limit $\eps \rightarrow 0$ to the system
\be \label{NRNR}
\left\{
\begin{array}{l}
\displaystyle - \frac{\partial }{\partial z} \left(A_z \frac{\partial \phi'}{\partial z} \right)  =  0 \,,\quad \text{in} \quad  \Omega\,, \\[3mm]
\displaystyle \frac{\partial \phi'}{\partial z} = 0 \quad \text{on}\quad
\Omega_x \times \partial \Omega_z \,, \qquad
\displaystyle \phi' = 0 \quad \text{on} \quad \partial \Omega_x \times \Omega_z\,,\\[3mm]
 \displaystyle \bar{\phi'} = 0 \quad \text{in} \quad  \Omega_x \,,
\end{array}
\right.
\ee
which is uniquely solvable, with the solution $\phi' \equiv 0$. Inserting this solution in (\ref{eq:compute:mean}), we conclude that
the solution of the AP formulation converges for $\eps \rightarrow 0$
towards the mean value part $\bar{\phi}(x)$, computed thanks to the
Limit model
\begin{equation}
(L)\,\,
\left\{
\begin{array}{l}
\displaystyle- \frac{\partial}{\partial x} \left( \bar{A}_\perp \frac{\partial \bar{\phi}}{\partial x}\right)   = \bar{f}(x) \,,\quad \text{in} \quad \Omega_x\,, \\[3mm]
 \displaystyle \bar{\phi} = 0 \quad \text{on} \quad \partial \Omega_x \,.
\end{array}
\right.
\end{equation}

The AP reformulation (\ref{eq:compute:mean})-(\ref{sys:phi:prime}) is equivalent to the Singular Perturbation problem 
(\ref{eq:ellipti:original_bis}) and is therefore valid for all
$0<\varepsilon <1$. This new formulation guarantees that, working with
a finite precision arithmetic, the computed solution converges in the
limit $\eps \rightarrow 0$ towards the solution of the limit model
(\ref{eq:ellipti:S}). This is a huge difference with the original
Singular Perturbation model which degenerates into an ill-posed
problem. Thus, by using the AP-formulation, \modif{we expect the computation of the numerical solution to be accurate}, uniformly in $\eps$.\\
  For the detailed mathematical proofs, we refer to the next section.

\begin{remark} \label{rem1} The condition $\overline{\phi'} =0$ in
  (\ref{sys:phi:prime}) holds automatically for $\eps>0$, since the right-hand side has zero average along the $z$-coordinate. Indeed, let $\psi$ be the solution of
\be \label{neige}
\left\{
\begin{array}{l}
\displaystyle - \frac{\partial }{\partial z} \left(A_z \frac{\partial \psi}{\partial z} \right) - \varepsilon \frac{\partial}{\partial x}\left(A_\perp \frac{\partial \psi}{\partial x}  \right) + \varepsilon \frac{\partial}{\partial x}\left(\overline{A_{\perp}' \frac{\partial \psi}{\partial x} }\right) = \varepsilon g' \,,\quad \text{in} \quad  \Omega\,, \\[3mm]
\displaystyle \frac{\partial \psi}{\partial z} = 0 \quad \text{on}\quad
\Omega_x \times \partial \Omega_z \,, \qquad
\displaystyle \psi = 0 \quad \text{on} \quad \partial \Omega_x \times \Omega_z\,,
\end{array}
\right.
\ee
with $\overline{g'}=0$. Taking the average along $z$, we get
$$
\left\{
\begin{array}{l}
\displaystyle- \frac{\partial}{\partial x} \left( \bar{A}_\perp \frac{\partial \bar{\psi}}{\partial x}\right)   = 0 \,,\quad \text{in} \quad \Omega_x\,, \\[3mm]
 \displaystyle \bar{\psi} = 0 \quad \text{on} \quad \partial \Omega_x \,,
\end{array}
\right.
$$
and thus $\overline{\psi} \equiv 0$, which is nothing but the
  constraint added in (\ref{sys:phi:prime}).\\
The computations of the
fluctuating part $\phi'$ via the equation \eqref{sys:phi:prime}
requires the discretization of an integro-differential operator. This
means that the discretization matrix will contain dense blocks.  However, using \eqref{eq:compute:mean} the system (AP2) can be rewritten as
\begin{equation}\label{sys:phi:prime:prime}
(AP2')\,\,\,
\left\{
\begin{array}{l}
\displaystyle - \frac{\partial }{\partial z} \left(A_z \frac{\partial \phi'}{\partial z} \right) - \varepsilon \frac{\partial}{\partial x}\left(A_\perp \frac{\partial \phi'}{\partial x}  \right)  = \\[3mm]
 \displaystyle \hspace{5.2cm} \varepsilon f + \varepsilon \frac{\partial}{\partial x} \left(A_\perp \frac{\partial \overline{\phi}}{\partial x} \right)\,,\quad \text{in} \quad  \Omega\,, \\[3mm]
 \displaystyle \frac{\partial \phi'}{\partial z} = 0 \quad \text{on}\quad
\Omega_x \times \partial \Omega_z \,, \qquad
\displaystyle \phi' = 0 \quad \text{on} \quad \partial \Omega_x \times \Omega_z\,, \\[3mm]
 \displaystyle  \overline{ \phi'}= 0\,, \quad \text{in} \quad \Omega_x\,.
\end{array}
\right.
\end{equation}
In this expression the right-hand side has no longer zero mean value
along the $z$-coordinate, but the integro-differential operator has
disappeared. The associated discretization matrix is thus sparser than
that obtained from the system~\eqref{eq:compute:mean}. Systems
(\ref{eq:compute:mean})-(\ref{sys:phi:prime}) and (\ref{eq:compute:mean})-(\ref{sys:phi:prime:prime})
are equivalent.
\end{remark}

\subsection{Mathematical study of the AP-formulation}\label{SEC32}
We establish in this section the mathematical framework of the
AP-formulation (\ref{eq:compute:mean})-(\ref{sys:phi:prime}) and study its
mathematical properties.
Let us thus introduce the two Hilbert-spaces
$$
\mathcal{V}:= \{ \psi(\cdot,\cdot) \in H^1(\Omega)\,\, / \,\, \psi=0 \,\,
\textrm{on} \,\, \partial \Omega_x \times \Omega_z \}\,, \quad \mathcal{W}:= \{ 
\psi(\cdot) \in H^1(\Omega_x)\,\, / \,\,  \psi=0 \,\,
\textrm{on} \,\, \partial \Omega_x \}\,, 
$$
with the corresponding scalar-products
\begin{equation} \label{sc_VW}
(\phi,\psi)_{\mathcal{V}}:=\varepsilon (\partial_x
\phi, \partial_x
\psi)_{L^2}+(\partial_z \phi,\partial_z \psi)_{L^2}\,, \quad (\phi,\psi)_{\mathcal{W}}:=(\partial_x
\phi,\partial_x
\psi)_{L^2}\,,
\end{equation}
\modif{ and the induced norms $||\cdot||_{\cal V}$, respectively $||\cdot||_{\cal W}$.}
For simplicity reasons, we denote in the sequel the $L^2$
scalar-product simply by the bracket $(\cdot,\cdot)$. Defining the
following bilinear forms
\begin{equation} \label{biVF}
\begin{array}{ll}
\displaystyle    a_0 \left(\phi', \psi' \right) &:= \displaystyle \int_0^{L_z} \int_{0}^{L_x} A_z(x,z)
   \frac{\partial \phi'}{\partial z}(x,z)  \frac{\partial \psi'}{\partial z}(x,z)
   dx dz  \, , \\[3mm]
\displaystyle    a_1 \left(\phi', \psi' \right) &:=\displaystyle  \int_0^{L_z} \int_{0}^{L_x} A_\perp(x,z)
    \frac{\partial \phi'}{\partial x}(x,z) \frac{\partial \psi'}{\partial x}
(x,z) dx dz  \, ,\\[3mm]
\displaystyle    a_2 \left(\bar \phi ,\bar \psi \right) &:=\displaystyle  \int_{0}^{L_x} \bar A_\perp(x)
  \frac{\partial \bar  \phi}{\partial x}(x)  \frac{\partial \bar
\psi}{\partial x} (x)  dx \, , \\[3mm]
 \displaystyle   c\left(\phi',\bar\psi \right) &:=\displaystyle   \int_0^{L_z} \int_{0}^{L_x} A_{\perp}'(x,z)
   \frac{\partial   \phi'}{\partial x}(x,z) \frac{\partial \bar  \psi}{\partial x}(x)  dx dz  \, ,\\[3mm]
 \displaystyle  d(\phi',\psi') &:=\displaystyle  \frac{1}{{L_z}} \int_{0}^{L_x} \int_0^{L_z} \int_0^{L_z}
  A_{\perp}'(x,z) \frac{\partial   \phi'}{\partial x}(x,z) \frac{\partial  \psi'}{\partial x}(x,\zeta)\,  dz d\zeta dx\, ,\\[3mm]
b(\bar P,\psi') &:= \displaystyle \int_{0}^{L_x} \bar P(x) \int_0^{L_z} \psi'(x,z) dz
  dx\, ,\\[4mm]
a(\phi',\psi')&:=\displaystyle a_0\left(\phi',\psi'\right) + \varepsilon
  a_1\left(\phi',\psi'\right) 
   - \varepsilon d(\phi',\psi')\,,
\end{array}
\end{equation}
permits to rewrite the AP system
(\ref{eq:compute:mean})-(\ref{sys:phi:prime}) under the weak form
\begin{equation}\label{eq:var:P}
(AP)\,\, 
\left\{
\begin{array}{l}
 \displaystyle a_2 \left(\bar\phi,\bar\psi \right)  = (\bar f,
\bar\psi)-  \frac{1}{L_z} c\left(\phi',\bar\psi
  \right)\, , \qquad \forall \bar \psi \in \mathcal{W} \, ,\\[3mm]
 \displaystyle  a\left(\phi',\psi'\right) + b(\bar P,\psi')= \varepsilon (f',\psi')- \varepsilon  c\left(\psi',\overline{\phi} \right)\, , \quad \forall 
  \psi' \in \mathcal{V}\, , \\[3mm]
  \displaystyle b(\bar Q,\phi') = 0 \, , \qquad \forall \bar Q\in \mathcal{W} \, ,
\end{array}
\right.  
\end{equation}
where $\phi'(x,z) \in \mathcal{V}$, $\bar \phi(x) \in \mathcal{W}$ as
well as $\bar{P}(x)
\in \mathcal{W}$ are the unknowns and $\psi' \in \mathcal{V}$, $\bar
\psi \in \mathcal{W}$ and $ \bar{Q} \in \mathcal{W}$ the test
functions. It can be observed that the constraint $\bar{\phi'}=0$ was
introduced via the Lagrange multiplier $\bar{P}$. We will see in the
next theorem that the weak formulation (\ref{eq:var:P}) is equivalent
for $\eps>0$ to the system
\begin{numcases}
  \strut a_2 \left(\bar\phi,\bar\psi \right)  = (\bar f,
\bar\psi)-  \frac{1}{L_z} c\left(\phi',\bar\psi
  \right)\, , \qquad \forall \bar \psi \in \mathcal{W} \, ,\label{eq:compute:mean:var} \\[3mm]
  a \left(\phi',\psi'\right) = \varepsilon (f',\psi')- \varepsilon  c\left(\psi',\overline{\phi} \right)\, , \quad \forall 
  \psi' \in \mathcal{V}\, ,\label{sys:phi:prime:var} 
\end{numcases}
where the explicit constraint $\bar{\phi'}=0$ does not appear. Let us assume in the sequel\\

\noindent {\bf Hypothesis A} {\it Let the diffusion functions $A_\perp \in L^{\infty} (\Omega)$ and
$A_z \in L^{\infty} (\Omega)$ satisfy
$$
0<c_\perp\le A_{\perp}(x,z) \le M_\perp \,,\quad  0<c_z\le A_{z}(x,z) \le M_z\,,
\quad \textrm{f.a.a.}\,\, (x,z) \in \Omega\,,
$$
with some positive constants $c_\perp, c_z, M_\perp, M_z$. Let moreover $f \in L^2(\Omega)$.}\\

\noindent The next theorem analyzes the well-posedness of the
AP-formulation.
\vspace{0.2cm}

\begin{theorem} \label{thm_EX}
For every $\eps>0$ the problem
(\ref{eq:compute:mean:var})-(\ref{sys:phi:prime:var}) admits under Hypothesis
A a unique solution
$(\phi'_\eps,\overline{\phi}_\eps) \in \mathcal{V} \times
\mathcal{W}$, where $\phi_\eps:=\phi'_\eps+\overline{\phi}_\eps$ is
the unique solution of the Singular Perturbation model
(\ref{eq:ellipti:original_bis}). The function $\phi'_\eps$ has zero
mean value along the $z$-coordinate, i.e. $\overline{\phi'}_\eps=0$
for every $\eps>0$.\\
Consequently, $(\phi'_\eps,\overline{\phi}_\eps)\in \mathcal{V} \times
\mathcal{W}$ is the unique solution of
(\ref{eq:compute:mean:var})-(\ref{sys:phi:prime:var}) if and only if $(\phi'_\eps,\overline{\phi}_\eps,\overline{P}_\eps)\in \mathcal{V} \times
\mathcal{W}\times \mathcal{W}$ is a solution of the AP-formulation (\ref{eq:var:P}). In this last case, we have $\overline{P}_\eps =0$.\\
Finally, these solutions satisfy the bounds
$$
||\phi_\eps||_{H^1(\Omega)} \le C ||f||_{L^2(\Omega)}\,, \quad ||\phi'_\eps||_{H^1(\Omega)}
\le C||f||_{L^2(\Omega)}\,, \quad ||\overline{\phi}_\eps||_{H^1(\Omega_x)} \le C||f||_{L^2(\Omega)}\,,
$$
with an $\eps$-independent constant $C>0$. In the limit $\eps \rightarrow 0$
there exist some functions $(\phi'_0,\overline{\phi}_0) \in \mathcal{V} \times
\mathcal{W}$, such that \modif{we have the following weak convergences in $H^1$}
\modif{$$
\phi'_\eps \rightharpoonup_{\eps \rightarrow 0} \phi'_0 \quad \textrm{in} \quad
H^1(\Omega)\,, \quad \overline{\phi}_\eps \rightharpoonup_{\eps \rightarrow 0}
\overline{\phi}_0\quad \textrm{in} \quad
H^1(\Omega_x)\,,
$$}
\modif{and the strong $L^2$ convergences}
\modif{$$
\phi'_\eps \rightarrow_{\eps \rightarrow 0} \phi'_0 \quad \textrm{in} \quad
L^2(\Omega)\,, \quad \partial_z \phi'_\eps \rightarrow_{\eps \rightarrow 0} \partial_z \phi'_0 \quad \textrm{in} \quad
L^2(\Omega)\,, \quad  \overline{\phi}_\eps \rightarrow_{\eps \rightarrow 0}
\overline{\phi}_0\quad \textrm{in} \quad
L^2(\Omega_x)\,,
$$}
where $\phi'_0 \equiv 0$ and $\overline{\phi}_0$ is the
unique solution of the Limit model (\ref{eq:ellipti:S}).
\end{theorem}

\noindent \debproof
The Singular Perturbation model (\ref{eq:ellipti:original_bis}) and the Limit model (\ref{eq:ellipti:S}) are standard
elliptic problems and posses under Hypothesis A (and for every $\eps >0$) unique solutions
$\phi_\eps \in \mathcal{V}$,  respectively
$\overline{\phi} \in \mathcal{W}$. It is then a simple consequence of the
decomposition (\ref{eq:def:decomp:bis}), that the problem
(\ref{eq:compute:mean:var})-(\ref{sys:phi:prime:var}) admits a
unique solution $(\phi'_\eps,\overline{\phi}_\eps) \in \mathcal{V} \times
\mathcal{W}$, where $\overline{\phi}_\eps(x):= \frac{1}{L_z}
\int_{0}^{L_z} \phi_\eps (x,z) dz$ is the mean and $\phi'_\eps:=\phi_\eps -\overline{\phi}_\eps$ the
fluctuation part. Thus we have also $\overline{\phi'_\eps}=0$. 
This property can also be understood from the
fact that
the right-hand side of (\ref{sys:phi:prime}), denoted in the sequel by
$g$
\begin{equation*}
  g(x,z):=f'(x,z)+\frac{\partial}{\partial x}
\left(A_\perp'(x,z) \frac{\partial \overline{\phi}}{\partial x}(x) \right)\,,
\end{equation*}
has zero mean value along the $z$-coordinate. Indeed, taking in 
(\ref{sys:phi:prime:var}) test functions $\psi'(x)
\in \mathcal{V}$ depending only on $x$, yields
immediately that $\overline{\phi'_\eps}=0$ for all $\eps >0$.\\
Standard
stability results for elliptic problems yield now the
$\eps$-independent estimate for the solution of the Singular Perturbation model
(\ref{eq:ellipti:original_bis})
$$
||\phi_\eps||^2_{H^1(\Omega)} \le ||\partial_x
\phi_\eps||^2_{L^2(\Omega)} + {1 \over \eps} ||\partial_z
\phi_\eps||^2_{L^2(\Omega)} \le C||f ||^2_{L^2(\Omega)} \,,
$$
implying that $||\overline{\phi}_\eps||^2_{H^1(\Omega_x)} \le C||f ||^2_{L^2(\Omega)}$ and
$||\phi'_\eps||^2_{H^1(\Omega)} \le C||f ||^2_{L^2(\Omega)}$, with a constant
$C>0$ independent of $\eps>0$. Thus there exist some functions $(\phi'_0,\overline{\phi}_0) \in \mathcal{V} \times
\mathcal{W}$, such that, up to a subsequence $\phi'_\eps
\rightharpoonup_{\eps \rightarrow 0} \phi'_0$ in $H^1(\Omega)$ and $\overline{\phi}_\eps \rightharpoonup_{\eps \rightarrow 0}
\overline{\phi}_0$ in $H^1(\Omega_x)$. Hence we have 
$$
\int_0^{L_x}\int_0^{L_z} \phi'_\eps(x,z) \psi(x,z) dx \, dz
\rightarrow_{\eps \rightarrow 0} \int_0^{L_x}\int_0^{L_z} \phi'_0(x,z)
\psi(x,z) dx \, dz\,, \quad \forall \psi \in \mathcal{V}\,.
$$
Taking here $\psi(x)\in \mathcal{V}$ depending only on the
$x$-coordinate, we observe that the feature $\overline{\phi'}_\eps \equiv 0$
yields the crucial property of the limit solution $\overline{\phi'}_0 \equiv 0$. Passing now to the limit $\eps
\rightarrow 0$ in (\ref{sys:phi:prime:var}),
we get that $\phi'_0$ is solution of
$$
\begin{array}{l}
 \displaystyle a_0(\phi_0',\psi')=0\,, \quad \forall \psi' \in \mathcal{V}\,,\quad \textrm{with} \quad 
\displaystyle  \overline{ \phi'_0}= 0\quad \text{in} \quad \Omega_x\,,
\end{array}
$$
which is the weak form of (\ref{NRNR}) and implies $\phi'_0 \equiv 0$. Finally, passing to the limit in (\ref{eq:compute:mean:var}), yields that $\overline{\phi}_0$ is the unique solution of
the Limit model (\ref{eq:ellipti:S}). Because of the uniqueness of the limit
$(\phi'_0,\overline{\phi}_0)$, we deduce that the whole sequence
$(\phi'_\eps,\overline{\phi}_\eps)$ converges weakly towards this
limit. To conclude the first part of the proof, we shall show \modif{the strong $L^2$ convergences}. For this, taking in
(\ref{sys:phi:prime:var}) $\phi_\eps'$ as test function and passing to
the limit $\eps \rightarrow 0$, yields $\partial_z \phi_\eps'
\rightarrow 0$ in $L^2(\Omega)$. As $\phi_\eps' \in \mathcal{V}$ and
$\bar{\phi_\eps'}=0$, the Poincar\'e inequality 
\modif{$$
||\phi_\eps'||_{L^2} \le C ||\partial_z \phi_\eps'||_{L^2}\,,
$$
is valid and implies that $\phi_\eps' \rightarrow 0$ in $L^2(\Omega)$. The convergence $\bar{\phi_\eps}\rightarrow \bar{\phi_0}$ in $L^2(\Omega_x)$ is immediate by compacity.}
It remains finally to prove
the equivalence between (\ref{eq:var:P}) and
(\ref{eq:compute:mean:var})-(\ref{sys:phi:prime:var}). This is
immediate. Indeed, if
$(\phi'_\eps,\overline{\phi}_\eps) \in \mathcal{V} \times \mathcal{W}$ is solution of
(\ref{eq:compute:mean:var})-(\ref{sys:phi:prime:var}), then
$(\phi'_\eps,\overline{\phi}_\eps,0)$ is solution of (\ref{eq:var:P}). And
if $(\phi'_\eps,\overline{\phi}_\eps,\overline{P}_\eps)\in \mathcal{V} \times \mathcal{W} \times \mathcal{W}$ satisfies
(\ref{eq:var:P}), then $\overline{P}_\eps \equiv 0$ (obvious by taking
as test function in (\ref{eq:var:P}) $\psi'(x) \in \mathcal{V}$ depending only on $x$) and $(\phi'_\eps,\overline{\phi}_\eps)$ solves hence
(\ref{eq:compute:mean:var})-(\ref{sys:phi:prime:var}).
\finproof

\noindent The subject of the next section will be the numerical resolution of
the AP-formulation  (\ref{eq:compute:mean})-(\ref{sys:phi:prime}) (or (\ref{eq:var:P})) and this
shall be done iteratively via a fixed-point application. Let us thus
introduce here the fixed-point map, construct an iterative
sequence and analyze its convergence. In the rest of this section, the
parameter $\eps>0$ shall be considered as fixed. Due to the fact that the two
systems (\ref{eq:var:P}) and (\ref{eq:compute:mean:var})-(\ref{sys:phi:prime:var})
are equivalent, we shall concentrate on the simpler one, i.e. (\ref{eq:compute:mean:var})-(\ref{sys:phi:prime:var}). Let us define the Hilbert space 
$$
\mathcal{U}:= \{ \psi(\cdot,\cdot) \in \mathcal{V}\,\, / \,\, \overline{\psi}=0\}\,,
$$
associated with the \modif{scalar product}
\modif{$$
(\phi,\psi)_*:=\int_0^{L_x}\!\!\int_0^{L_z}\!\!{A}_z \partial_z \phi
\, \partial_z \psi
dz dx+ \eps \int_0^{L_x}\!\!\int_0^{L_z}\!\!{A}_\perp \partial_x
\phi\, \partial_x \psi dz
dx\,,
$$}
\modif{which is equivalent to the scalar product $(\cdot,\cdot)_{\cal V}$ on
$\cal V$, defined by (\ref{sc_VW}).}\\
The fixed-point map $T: \mathcal{U} \rightarrow \mathcal{U}$
is defined as follows: With $\phi' \in \mathcal{U}$ we associate $\overline{\phi} \in \mathcal{W}$,
solution of (\ref{eq:compute:mean:var}). Then constructing the right-hand side
of (\ref{sys:phi:prime:var}) via this $\overline{\phi} \in \mathcal{W}$, we
define $T(\phi')$ as the corresponding solution of
(\ref{sys:phi:prime:var}). Denoting by $(\phi'_*,\overline{\phi}_*) \in
\mathcal{V} \times \mathcal{W}$ the unique solution of
(\ref{eq:compute:mean:var})-(\ref{sys:phi:prime:var}), we remark by
Theorem \ref{thm_EX} that $\phi'_* \in \mathcal{U}$ and that it is the unique fixed-point of the
map $T$.\\

\begin{theorem} \label{FP_thm}
Let $\eps>0$ be fixed and let $\phi'_* \in \mathcal{U}$ be the unique fixed-point of the
application $T: \mathcal{U} \rightarrow \mathcal{U}$ constructed as
follows
$$
\phi'\in \mathcal{U} \quad \xrightarrow{(\ref{eq:compute:mean:var})}
\quad \overline{\phi}\in \mathcal{W}
\quad \xrightarrow{(\ref{sys:phi:prime:var})} \quad  T(\phi')\in \mathcal{U}\,.
$$
Then for every starting point $\phi'_0 \in \mathcal{U}$,  the sequence $\phi'_k:=T(\phi'_{k-1}) =
T^k(\phi'_0)$ converges in \modif{$(\mathcal{U},||\cdot||_*)$, and
  consequently also in $(\mathcal{U},||\cdot||_{\cal V})$,} towards the fixed-point
$\phi'_* \in  \mathcal{U}$ of $T$.
\end{theorem}

\noindent The proof of this theorem is based on the following\\

\begin{lemme} \cite{brezis} \label{BR}
Let $(\mathcal{U},||\cdot||_{\modif{*}})$ be a normed space and $T: \mathcal{U} \rightarrow
\mathcal{U}$ a contractive application, i.e.
$$
||T(\phi)-T(\psi)||_{\modif{*}} <||\phi-\psi||_{\modif{*}}\,,
\quad \forall \phi, \psi \in  \mathcal{U}\quad \textrm{with} \quad \phi \neq \psi\,.
$$
Then the set of fixed-points of T, denoted by $FP(T)$, is identical
with the set of accumulation points of the sequences $\{
T^k(\phi)\}_{k \in \NN}$,
with $\phi \in \mathcal{U}$, set which is denoted by $AP(T)$. Moreover, these two spaces contain at most one
element.
\end{lemme}

\noindent {\bf Proof of theorem \ref{FP_thm} :}\\ 
The linear application $T$ is well-defined. The first step $\phi'\in \mathcal{U} \quad \xrightarrow{(\ref{eq:compute:mean:var})}
\quad \overline{\phi}\in \mathcal{W}$ is immediate by the Lax-Milgram theorem. For the second step, we remark that for given $\overline{\phi} \in \mathcal{W}$ the equation
\be \label{AP_c}
a(\theta,\psi')=\varepsilon (f',\psi')-\varepsilon c(\psi',\overline{\phi})\,, \quad \forall \psi' \in \mathcal{V}\,,
\ee
has a unique solution $\theta \in \mathcal{U}$. Indeed, we notice first (by taking test functions only depending on the $x$-coordinate) that $\overline{\theta}=0$. This enables us to consider instead of (\ref{AP_c}), the variational
formulation
\be \label{E_APm}
m(\theta,\psi')=\varepsilon (f',\psi')-\varepsilon c(\psi',\overline{\phi})\,, \quad \forall \psi' \in \mathcal{V}\,,
\ee
where the bilinear form $a(\cdot,\cdot)$, which is not coercive, was replaced by the coercive bilinear form $m(\cdot,\cdot)$, given by
\be \label{E_APmm}
m(\theta,\psi'):=a(\theta,\psi')+{ \varepsilon M_\perp \over L_z} \int_0^{L_x}
\left[ \int_0^{L_z} \partial_x \theta (x,z) dz \right]\left[
\int_0^{L_z} \partial_x \psi' (x,z) dz \right]  dx \,.
\ee
Indeed, due to the
property $\overline{\theta}=0$, the two equations (\ref{AP_c}) and
(\ref{E_APm}) are equivalent and this time $m(\cdot,\cdot)$ is a continuous,
coercive bilinear form, as for all $\psi' \in \mathcal{V}$
we have
$$
m(\psi',\psi')\ge \int_0^{L_x} \int_0^{L_z} A_z |\partial_z \psi'|^2\, dz
dx + \varepsilon \int_0^{L_x} \int_0^{L_z} A_\perp |\partial_x
\psi'|^2\, dz dx\ge C||\psi'||_{\mathcal{V}}^2\,.
$$
Thus the Lax-Milgram theorem implies the existence and uniqueness of a
solution $\theta \in \mathcal{U}$ of the continuous problem
(\ref{E_APm}) and hence also of
problem (\ref{AP_c}).  We have shown by this that $T$ is a well-defined mapping.\\
Furthermore we know that $T$ admits, for fixed $\eps>0$, a unique fixed-point, denoted by $\phi'_* \in \mathcal{U}$. Let us now suppose that we have shown that $T$
is contractive. Then lemma \ref{BR} implies that
$FP(T)=AP(T)=\{\phi'_*\}$. Thus choosing an arbitrary starting point
$\phi'_0 \in \mathcal{U}$, and
constructing the sequence $\phi'_k:=T(\phi'_{k-1}) =
T^k(\phi'_0)$, we deduce that this sequence has a unique accumulation
point $\phi'_*$ in $\mathcal{U}$. This means that the sequence $\{\phi'_{k}\}_{k\in
\NN}$ converges in \modif{$(\mathcal{U},||\cdot||_*)$ towards
$\phi'_*$. Due to the fact that $||\cdot||_*$ and $||\cdot||_{\cal V}$
are equivalent norms, we have also the convergence in $(\mathcal{U},||\cdot||_\mathcal{V})$.}\\

\noindent It remains to show that $T$ is contractive. For this let $\phi'_1,
\phi'_2 \in \mathcal{U}$ be two given, distinct functions. Denoting by
$\phi':=\phi'_1-\phi'_2$,
$\overline{\phi}:=\overline{\phi}_1-\overline{\phi}_2$ (where
$\overline{\phi_i} \in \mathcal{W}$ are the corresponding solutions of (\ref{eq:compute:mean:var})) and $\theta':=
T(\phi'_1)-T(\phi'_2)$, we have to show that $||\theta'||_{\modif{*}} <
||\phi'||_{\modif{*}}$. First we observe that $\overline{\phi}$ solves
\be \label{EQ_T1}
a_2(\bar \phi,\bar \psi)=-{1 \over L_z} c(\phi',\bar \psi)\,, \quad \forall \bar \psi \in \mathcal{W}\,,
\ee
and $\theta'$ is solution of
\be \label{EQ_T2}
a(\theta',\psi')=- \varepsilon  c\left(\psi',\overline{\phi} \right)\, , \quad \forall 
  \psi' \in \mathcal{V}\, .
\ee
Taking in (\ref{EQ_T1})
$\overline{\phi}$ as test function, gives rise to
$$
\begin{array}{lll}
\displaystyle \int_0^{L_x} \overline{A}_\perp |\partial_x \overline{\phi} (x)|^2 \,
dx &=& \displaystyle - \int_0^{L_x} \left[ {1 \over L_z} \int_0^{L_z}
{A}'_\perp \partial_x \phi' (x,z) dz \right] \partial_x
\overline{\phi} (x)\, dx\\[5mm]
&= &\displaystyle - \int_0^{L_x} \left[ {1 \over L_z} \int_0^{L_z}
{A}_\perp \partial_x \phi' (x,z) dz \right] \partial_x
\overline{\phi} (x)\, dx\\[5mm]
&\le &\displaystyle {1 \over \sqrt{L_z}} \left[ \int_0^{L_x}\!\! \int_0^{L_z}\!\!{A}_\perp |\partial_x \phi'|^2 dz dx
\right]^{1/2} \left[ \int_0^{L_x}\!\!\overline{A}_\perp |\partial_x \overline{\phi}|^2 dx
\right]^{1/2} \,.
\end{array}
$$
Thus
$$
\left[ \int_0^{L_x} \overline{A}_\perp |\partial_x \overline{\phi}
(x)|^2 dx \right]^{1/2} \le {1 \over \sqrt{L_z}} \left[
\int_0^{L_x}\int_0^{L_z}{A}_\perp |\partial_x \phi'|^2 dz dx
\right]^{1/2}\,.
$$
Equally, taking in (\ref{EQ_T2}) $\theta'$ as test function gives rise
to

\begin{equation}\label{NR}
  \begin{split}
   \int_0^{L_x}\!\!\int_0^{L_z}\!\!&\!\!{A}_z |\partial_z \theta'|^2 dz dx
   + \eps \int_0^{L_x}\!\!\int_0^{L_z}\!\!\!\!{A}_\perp |\partial_x
   \theta'|^2 dz  
    \le - \eps \int_0^{L_x}\!\!\int_0^{L_z}\!\!\!\!{A}_\perp \partial_x
   \overline{\phi}\, \partial_x \theta'  dz dx \hspace*{-2cm}\\[5mm]
   & \le  \eps \left[ \int_0^{L_x}\!\!\int_0^{L_z}\!\!{A}_\perp
     |\partial_x \overline{\phi}|^2 dz dx \right]^{1/2} \!\!\left[ 
     \int_0^{L_x}\!\!\int_0^{L_z}\!\!{A}_\perp |\partial_x \theta'|^2
     dz dx\right]^{1/2}\\[5mm]
   &\le \displaystyle \eps \sqrt{L_z} \left[ \int_0^{L_x}\!\!\overline{A}_\perp |\partial_x \overline{\phi}|^2 
     dx \right]^{1/2} \!\!\left[ \int_0^{L_x}\!\!\int_0^{L_z}\!\!{A}_\perp |\partial_x \theta'|^2 dz
     dx\right]^{1/2}\,.
  \end{split}
\end{equation}
This last inequality yields
$$
\begin{array}{lll}
\displaystyle \int_0^{L_x}\!\!\int_0^{L_z}\!\!{A}_z |\partial_z \theta'|^2
dz dx&+& \displaystyle \eps \int_0^{L_x}\!\!\int_0^{L_z}\!\!{A}_\perp |\partial_x \theta'|^2 dz
dx \le \eps L_z \int_0^{L_x}\!\!\overline{A}_\perp |\partial_x \overline{\phi}|^2 
dx\\[3mm]
&\le & \displaystyle \eps \int_0^{L_x}\!\!\int_0^{L_z}\!\!{A}_\perp |\partial_x \phi'|^2 dz
dx\\[3mm]
&< & \displaystyle \int_0^{L_x}\!\!\int_0^{L_z}\!\!{A}_z |\partial_z \phi'|^2
dz dx+ \eps \int_0^{L_x}\!\!\int_0^{L_z}\!\!{A}_\perp |\partial_x \phi'|^2 dz
dx\,.
\end{array}
$$
In this last step we would have the ``equality'' if and only if $\int_0^{L_x}\int_0^{L_z}{A}_z |\partial_z \phi'|^2
dz dx =0$. This is however only possible for functions depending exclusively
on the $x$-coordinate, $\phi'(x)$, which is in contradiction with the
fact that $\overline{\phi'}=0$ and $\phi'\neq 0$. Thus we have shown that $||T(\phi')||_{\modif{*}} <
||\phi'||_{\modif{*}}$ for $\phi' \neq 0$, $\phi' \in \mathcal{U}$,
which means that $T$ is a contractive application on \modif{$(\mathcal{U},||\cdot||_*)$.}
\finproof
\section{Numerical discretization and simulation results}\label{SEC3}
This part of the paper is concerned with the numerical
discretization of the AP-scheme
(\ref{eq:compute:mean})-(\ref{sys:phi:prime}) and the
comparison of the simulation results with those obtained via the
Singular Perturbation model (\ref{eq:ellipti:original_bis}) and the Limit model (\ref{eq:ellipti:S}).
\subsection{Discretization}\label{SEC31}
The numerical resolution of the Asymptotic Preserving system
(\ref{eq:compute:mean})-(\ref{sys:phi:prime}) is done by means
of the standard finite element method.\\
Let us recall the variational formulation of the AP-formulation
\begin{equation}\label{eq:var:P:bis}
\left\{
\begin{array}{l}
 \displaystyle a_2 \left(\bar\phi,\bar\psi \right)  = (\bar f,
\bar\psi)-  \frac{1}{L_z} c\left(\phi',\bar\psi
  \right)\, , \qquad \forall \bar \psi \in \mathcal{W} \, ,\\[3mm]
 \displaystyle  a\left(\phi',\psi'\right) + b(\bar P,\psi')= \varepsilon (f',\psi')- \varepsilon  c\left(\psi',\overline{\phi} \right)\, , \quad \forall 
  \psi' \in \mathcal{V}\, , \\[3mm]
  \displaystyle b(\bar Q,\phi') = 0 \, , \qquad \forall \bar Q\in \mathcal{W} \, ,
\end{array}
\right.  
\end{equation}
with the notation of section \ref{SEC2}.
Here $\phi'(x,z) \in \mathcal{V}$, $\bar \phi(x) \in \mathcal{W}$ as
well as $\bar{P}(x)
\in \mathcal{W}$ are the unknowns and $\psi' \in \mathcal{V}$, $\bar
\psi \in \mathcal{W}$ and $ \bar{Q} \in \mathcal{W}$ the test functions. \\
The introduction of the Lagrange multiplier $\overline{P}(x)$ was
explained in a simplistic manner in the preceding sections and will be analyzed in more details in section \ref{SEC4}. Due to the equivalence of
(\ref{eq:var:P:bis}) and (\ref{eq:compute:mean:var})-(\ref{sys:phi:prime:var}), one can comment that the introduction of $\overline{P}(x)$ is superfluous, but this is not the case for the discretized equations. The property $\overline{\phi'}=0$ is indeed automatically fulfilled since the right-hand side of equation \eqref{sys:phi:prime:var} has a zero mean value along the $z$-coordinate. However the discrete implementation of this quantity introduces round-off errors which probably will destroy the zero mean value property and justify the introduction of the Lagrange multiplier. \\
For simplicity reasons we
omitted here the $\eps$-index of the solution
$(\phi'_\eps,\overline{\phi}_\eps)$, the parameter $\eps >0$ being considered as fixed.\\

\noindent To discretize now the
system (\ref{eq:var:P:bis}) we introduce the grid 
$$
0=x_0 \le \cdots \le x_n \le \cdots \le x_{N_x+1}=L_x\,, \quad 0=z_1 \le
\cdots \le z_k \le \cdots \le z_{N_z}=L_z
$$
and denote the cells by
$I_n:=[x_{n},x_{n+1}]$ and $J_k:=[z_{k},z_{k+1}]$. The finite
dimensional spaces $\mathcal{V}_h \subset \mathcal{V}$ and $\mathcal{W}_h \subset \mathcal{W}$ are constructed
as usual, by
means of the hat functions ($\mathcal{Q}_1$ finite elements)
$$
\chi_n(x):= \left\{
\begin{array}{ll}
\displaystyle \frac{x-x_{n-1}}{x_{n}-x_{n-1}}\,,& x \in I_{n-1}\,,\\[2mm]
\displaystyle \frac{x_{n+1}-x}{x_{n+1}-x_{n}}\,,& x \in I_n\,,\\[2mm]
\displaystyle 0 \,,& \textrm{else}
\end{array} \right. \,, \quad 
\kappa_k(x):= \left\{
\begin{array}{ll}
\displaystyle \frac{z-z_{k-1}}{z_{k}-z_{k-1}}\,,& z \in J_{k-1}\,,\\[2mm]
\displaystyle \frac{z_{k+1}-z}{z_{k+1}-z_{k}}\,,& z \in J_k\,,\\[2mm]
\displaystyle 0 \,,& \textrm{else}
\end{array} \right. \,.
$$
Thus we are searching for approximations $\phi'_h \in \mathcal{V}_h$,
$\bar \phi_h \in \mathcal{W}_h$ and $\bar P_h \in \mathcal{W}_h$, which can be written
under the form
$$
\phi'_h(x,z)=\sum_{n=1}^{N_x} \sum_{k=1}^{N_z} \alpha_{nk} \chi_n(x)
\kappa_k(z) \,, \quad \bar \phi_h(x)=\sum_{n=1}^{N_x}  \beta_{n} \chi_n(x) \,, \quad
\bar P_h(x)=\sum_{n=1}^{N_x}  \gamma_{n} \chi_n(x)\,.
$$
Inserting these decompositions in the variational formulation
(\ref{eq:var:P:bis}) and taking as test functions the hat-functions
$\chi_n$ and $\kappa_k$ gives rise to the following linear system to be
solved in order to get the unknown coefficients $\alpha_{nk}$, $\beta_n$ and
$\gamma_n$
\begin{align}\label{eq:Dyn:moy}
    A_2 \beta &= {\bf w} \,,\\
    \left(  \begin{array}[c]{cc}
    A_0 + \varepsilon \left(A_1-D\right)  &B \\
    B^t & 0 
  \end{array}\right) 
    \left( 
  \begin{array}[c]{c}
    \alpha \\ \gamma
  \end{array}
  \right) &=\varepsilon \left( 
  \begin{array}[c]{c}
    {\bf v}\\ 0
  \end{array}
  \right)\,,\label{eq:Dyn:fluc}
\end{align}
where the matrices $A_2 \in \RR^{N_x \times N_x}$, $A_0, A_1, D \in
\RR^{N_x N_z \times N_x N_z}$ and $B \in \RR^{N_x N_z \times N_x}$ correspond to the bilinear forms (\ref{biVF}) and the right-hand
sides are defined by
$$
{\bf w}_n:=(\bar f, \chi_n)-\frac{1}{L_z} c(\phi'_h,\chi_n)\,, \quad
{\bf v}_{nk}:=(f', \chi_n \kappa_k)-c(\chi_n \kappa_k,\bar \phi_h) = (g, \chi_n \kappa_k)\,,
$$
for all $\quad n=1,\cdots,N_x\,; \,\,\, k=1, \cdots N_z\,$ and
\begin{equation}\label{def:function:g}
   g(x,z):=f'(x,z)+\frac{\partial}{\partial x}
\left(A_\perp'(x,z) \frac{\partial \overline{\phi}}{\partial x}(x) \right)\,.
\end{equation}
Solving iteratively the linear systems (\ref{eq:Dyn:moy})-(\ref{eq:Dyn:fluc})
permits finally to get the unknown function $\phi_h(x,z)=\bar \phi_h (x) +
\phi'_h(x,z)$. The convergence of the iterations was proved for the
continuous case in theorem
\ref{FP_thm} and can be identically adapted for the discrete case.
\subsection{Numerical results}\label{SEC33}
In this section we shall compare the numerical results obtained by the
discretization of
the Singular Perturbation model, the Limit model and the just presented
Asymptotic Preserving reformulation. With this aim, we consider a test
case where the exact solution is known. Let thus
\begin{gather} \label{SOL_EX}
  \phi_e(x,z) := \sin\left(\frac{2\pi}{L_x} x\right)+ \varepsilon \cos\left(\frac{2\pi}{L_z} z\right)\sin\left(\frac{2\pi}{L_x} x\right) \,,
\end{gather}
be the exact solution of  problem (\ref{eq:ellipti:original_bis}),
where we choose $  A_\perp(x,y) = c_1 + x z^2$ and $A_z(x,z) = c_2 +
xz$,  with two constants $c_1>0$, $c_2>0$. \modif{The numerical
experiments are performed with $L_x=L_z= 10$ and $c_1 = c_2 = L_z$.} The exact right-hand side $f$ is computed by inserting (\ref{SOL_EX}) in (\ref{eq:ellipti:original_bis}).
We denote by $\phi_P$, $\phi_L$
and $\phi_A$, respectively, the numerical solutions of the Singular
Perturbation model \eqref{eq:ellipti:original_bis}, the Limit model
\eqref{eq:ellipti:S}  and the Asymptotic Preserving formulation (\ref{eq:compute:mean})-(\ref{sys:phi:prime}). The comparison will be
done in the $l^{2}$-norm, that means
\begin{equation}
  ||  \phi_e - \phi_{num}||_2 = {1 \over \sqrt{N}} \left( \sum_{i\in\mathcal{G}} \left| \phi_e(X_i) - \phi_{num,i} \right|^2 \right) ^{1/2} \,,
\end{equation}
where $\phi_{num}$ stands for one of the numerical solutions and
$\phi_e(X_i)$ is the exact solution evaluated in the grid point
$X_i$. The index $i$ covers all possible grid indices, reassembled in
the set $\mathcal{G}$, and $N$ is the total number of grid points.
The linear systems obtained after the discretization  of either the P-model, the L-model or the AP-formulation are solved thanks to the same numerical algorithm (MUMPS \cite{MUMPS}). The purpose here is not to design a specific preconditioner for the resolution of these linear systems, but to point out the efficiency of the presently introduced AP-method to deal with a large range of anisotropy ratios.

\begin{figure}[htpb]
  \centering
  \psfrag{Dynamo}[l][l][.5]{(S)}
  \psfrag{Striation}[l][l][.5]{(L)}
  \psfrag{Asymptotic Preserving}[l][l][.5]{(AP)}
  \psfrag{Epsilon}[][][.8]{$\varepsilon$}
  \psfrag{Error Norm 2}[][][1.]{}

  \subfigure[Grid with $50 \times 50$ points.\label{fig:error:coarse}]{
  \begin{minipage}[c]{0.45\textwidth}
    \includegraphics[width=\textwidth]{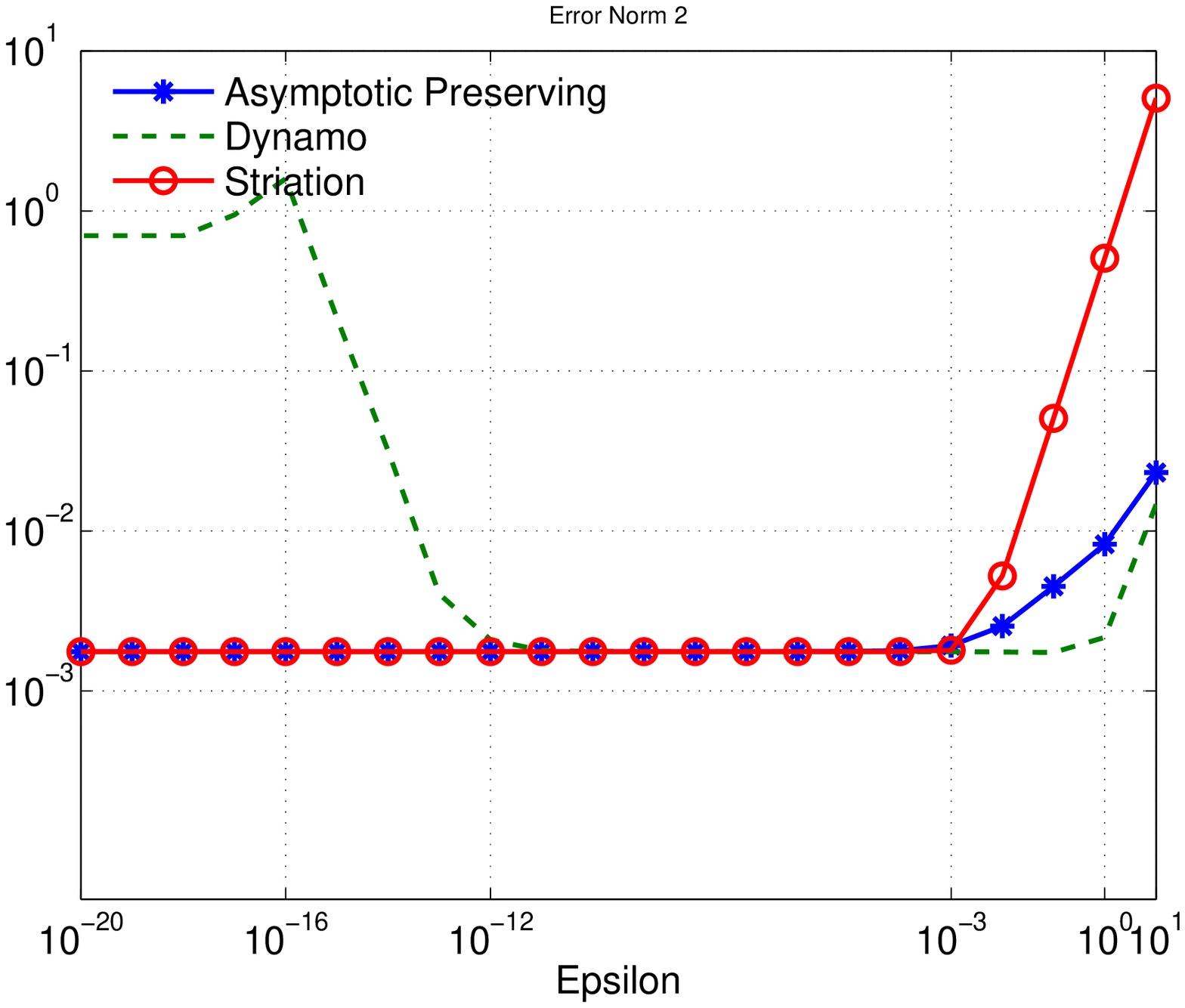}
    \abovecaptionskip 0cm
  \end{minipage} }\hfill%
  \subfigure[Grid with $500 \times 500$ points.\label{fig:error:refined}]{
  \begin{minipage}[c]{0.45\textwidth}
    \psfrag{Epsilon}[][][1.]{$\varepsilon$}
    \psfrag{Error Norm Infinity}[][][1.]{}
    \includegraphics[width=\textwidth]{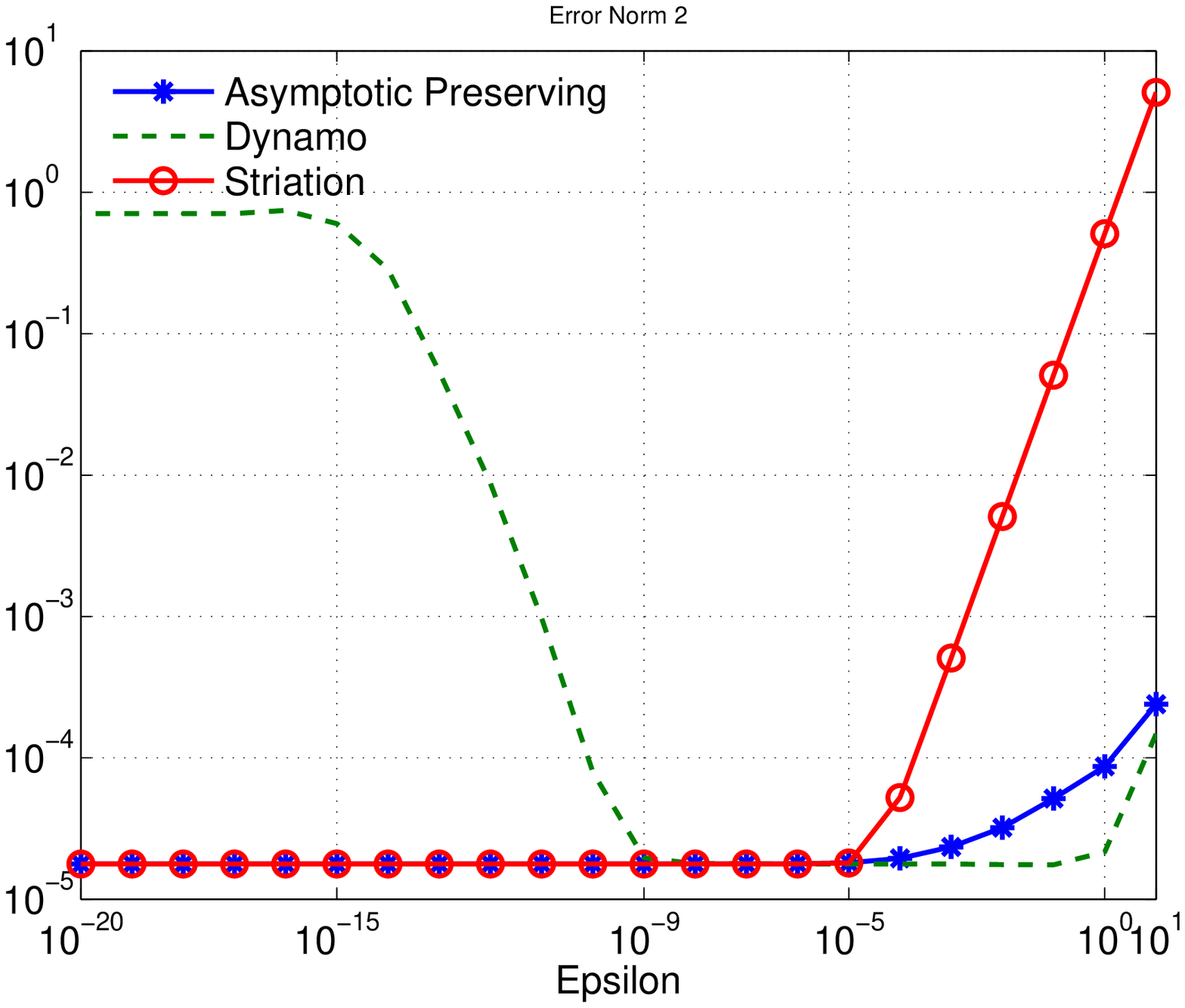}
  \end{minipage} }

\abovecaptionskip 0cm
  \caption{Absolute error in the $l^{2}$-norm between the
computed solutions $\phi_P, \phi_L, \phi_A$ and the exact solution
$\phi_e$, as a function of $\varepsilon$ and on different grids.
Dashed lines~: (S) Standard scheme~: discretization of the P-model;
Stars~: (AP) AP-scheme; Circles~: (L) discretization
of the L-model.}
  \label{fig:error:30}
\end{figure}

\noindent As can be seen from Table \ref{tab:error} and Figure
\ref{fig:error:30}, the finite element resolution of the Singular Perturbation model
is precise only for large $0<\eps<1$, whereas the Limit model is
accurate for small $\eps \ll 1$.
\begin{table}[!ht] 
  \centering
  \begin{tabular}[c]{|c||c|c|c|c|c|c|} \hline
    $\varepsilon$ & $10$   & $1$          & $10^{-1}$        & $10^{-4}$     & $10^{-14}$& $10^{-16}$ \\ \hline\hline
    AP-scheme & $3.4 \cdot  10^{-2}$  &$7.8 \cdot 10^{-3}$& $3.8 \cdot 10^{-3}$&$2.7 \cdot 10^{-3}$&$2.7 \cdot 10^{-3}$& $2.7 \cdot 10^{-3}$ \\ \hline
    S-scheme & $2.8\cdot10^{-2}$&$4.5 \cdot 10^{-3}$& $2.8 \cdot 10^{-3}$&$2.7 \cdot 10^{-3}$&$6.6 \cdot 10^{-2}$& $1.2 $ \\ \hline
    L-model & $9.9 $   &$1.0 \cdot 10^{1}$ & $1.0 \cdot 10^{-1}$&$2.8 \cdot 10^{-3}$&$2.7 \cdot 10^{-3}$& $2.7 \cdot 10^{-3}$ \\ \hline
  \end{tabular}
  \caption{Absolute error in the $l^{\infty}$-norm for the approximation computed thanks to the AP-scheme, discretized Singular Perturbation and Limit models (S-scheme and L-model) as compared to the exact solution.}
  \label{tab:error}
\end{table}
The range of $\eps$-values in which both the Singular Perturbation and
the Limit models provide an accurate approximation of the solution
shrinks as the mesh size is refined. For a coarse grid (with $50\times
50$ points see figure~\ref{fig:error:coarse}) this domain ranges from
$10^{-12}$ to $10^{-3}$ while it is reduced to $10^{-9} - 10^{-5}$ for
the refined $500 \times 500$ grid
(figure~\ref{fig:error:refined}). This question is determinant for the
development of a model coupling strategy. Indeed it requires an
  intermediate area where both discretized models furnish an accurate
  approximation and we observe that for refined meshes this area may
  not exist. This reduction of the validity domain can be explained for both the L-model and P-model but for quite different reasons.

The numerical approximation computed via the
Limit model is altered by both the discretization error of the
numerical scheme and the approximation error introduced by the
reduction of the initial Singular Perturbation problem to the Limit
problem. For coarse grids, the global error is rapidly dominated by
the scheme discretization error, but as the mesh is refined, the
approximation error  becomes preponderant, as the Limit model is
precise only for small $\eps$-values. The schemes implemented here are
of second order, thus when the mesh size is divided by ten, the 
discretization error is reduced by one hundred. The global error for
the L-model displayed in figure~\ref{fig:error:coarse} does not depend
on $\eps$ as soon as $\eps< 10^{-3}$. Below this limit the L-model is
able to furnish a better approximation of the solution with vanishing
$\eps$, however the numerical scheme is not precise enough and
consequently the global error does not decrease. For the refined mesh, this discretization error is lowered by two order of magnitudes and the global error is a function of $\eps$ as long as its value is greater than $ 10^{-5}$ (Fig.~\ref{fig:error:refined}). 

\begin{figure}[!ht]
  \centering
  \begin{minipage}[c]{0.45\textwidth}
    \centering
    \psfrag{Epsilon}[][][1.]{$\varepsilon$-values}
    \psfrag{Error}[][][1.]{Condition number estimate}

    \psfrag{Dynamo 1}[l][l][.7]{S-scheme (50$\times$50)}
    \psfrag{Dynamo 2}[l][l][.7]{S-scheme (200$\times$200)}
    \psfrag{Asymptotic Preserving 1}[l][l][.7]{AP-scheme (50$\times$50)}
    \psfrag{Asymptotic Preserving 2}[l][l][.7]{AP-scheme (200$\times$200)}  
    \includegraphics[width=\textwidth]{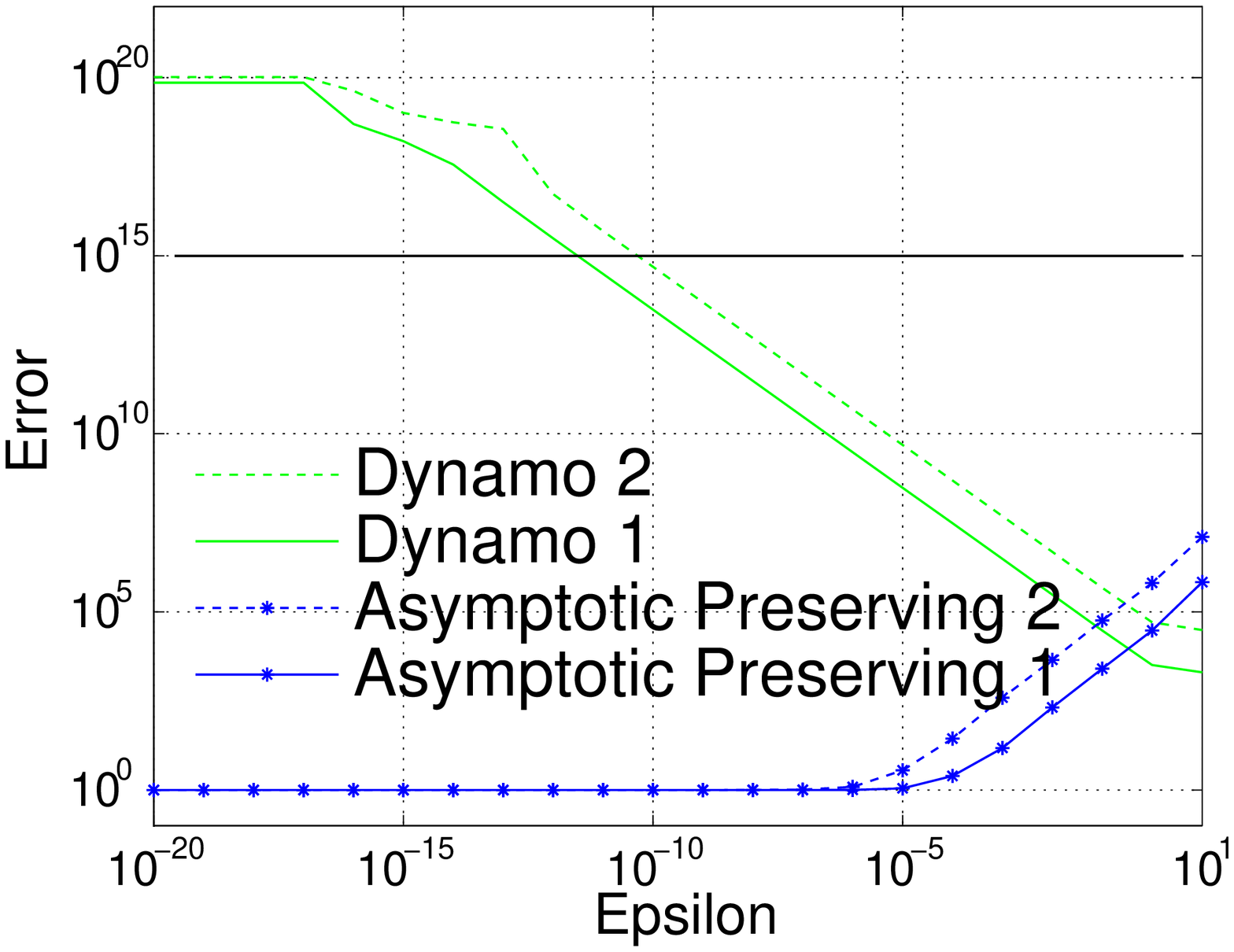}
    \caption{Condition number estimate for the discretization matrices of the Standard (S) and AP schemes (computed by LAPACK \cite{LAPACK}) as a function of $\eps$. Different grids of  50$\times$50 and 200$\times$200 points and different $\varepsilon$-values are used. Dashed/Plain lines~: $200\times 200$ / $50\times50$ grid ; Stars~: AP-scheme.}
    \label{fig:condition:number}    
  \end{minipage}\hfill%
  \begin{minipage}[c]{0.45\textwidth}
  \centering
  \psfrag{Fluctuation}[l][l][.8]{$\left\| \phi_A' -\phi_e' \right\|_2 $}
  \psfrag{Mean}[l][l][.8]{$\left\| \bar \phi_A -\bar \phi_e \right\|_2 $}
  \psfrag{Error}[][][1.]{$l^2$-Error}
  \psfrag{Iterations}[][][1.]{Iterations}
  \includegraphics[width=\textwidth]{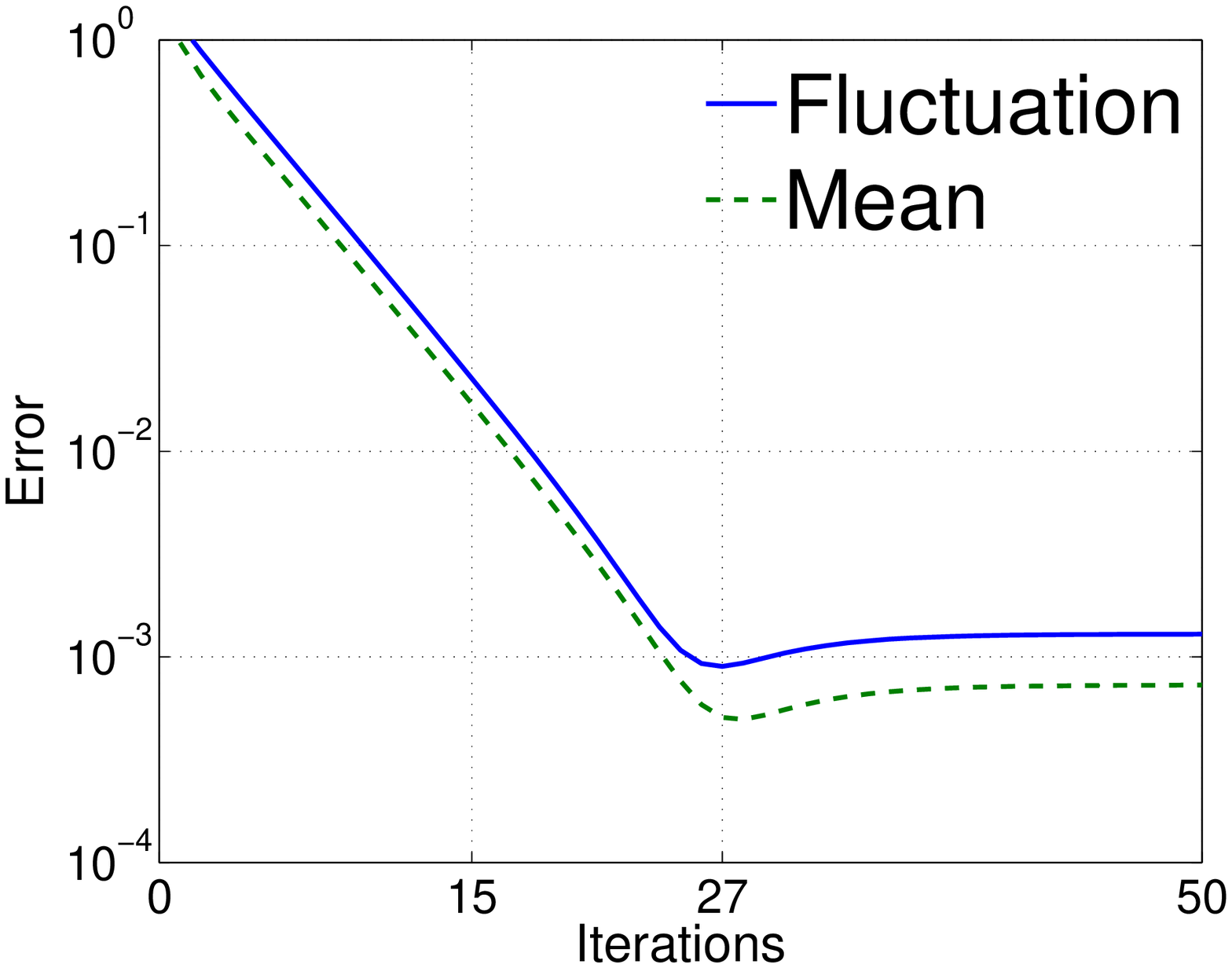}
  \caption{The $l^2$ absolute error between the exact solution and the numerical approximation computed with the AP-scheme, as a function of the iteration number, with $\eps=10$ and a $200 \times 200$-mesh. Dashed line~: mean part of the solution; Plain line~: fluctuating part.}
 \label{fig:relaxation}
  \end{minipage}
\end{figure}

The analysis for the Singular Perturbation model is quite complementary. The accuracy of the approximation provided by the P-model is good for large
$\eps$-values and deteriorates rapidly for small ones. This can be
explained by the conditioning of the linear system obtained by the
P-model discretization. An estimate of the condition number for the
matrix is displayed in figure~\ref{fig:condition:number} for two
different grid sizes. This conditioning deteriorates with vanishing
$\eps$-parameter, which is coherent with the fact that, working with a finite-precision arithmetic, the Singular Perturbation model degenerates into an ill-posed problem.
 This also explains
the blow up of the error displayed in figure~\ref{fig:error:30} as
soon as the conditioning of the matrix approaches the critical value
of the double precision (materialized by the level $10^{15}$ in
Fig.~\ref{fig:condition:number}). This limit is reached on more
refined meshes for larger $\eps$-values  ($\eps\approx10^{-12}$ on a
$50\times 50$ grid and $\eps \approx 10^{-10}$ on a $200\times200$
grid). As expected, the P-model, though valid for all $\eps$-values,
cannot be exploited numerically for small $\eps$. The $\eps$-region
where both the P-model and the L-model are accurate all-together,
shrinks dramatically with the size of the mesh, fact which motivates the development of the AP-method.\\
The condition number estimate of the linear system providing the
approximation of the solution for the AP-scheme is also plotted in
Figure~\ref{fig:condition:number}. The conditioning of the system is
rather $\eps$ independent and this is due to the introduction of the
Lagrange multiplier, which forces the system in the limit to remain
well-posed. The accuracy of the AP-scheme is totally comparable to the
P-model for the large values of $\eps$ and to the L-model for the
smallest ones. 
The AP-formulation is a good tool for computing an approximation for
the solution which is accurate uniformly in $0<\eps<1$ and is
therefore of great practical interest. Note that this approximation is
obtained thanks to an iterative sequence  $\{ \phi'_k \}_{k \in \NN}$,
constructed with the fixed-point mapping $T$ defined in theorem
\ref{FP_thm}. The convergence of this iterative process is analysed in
figure~\ref{fig:relaxation} on a $200\times 200$ grid for a large
value of $\eps$. The $l^2$-absolute error between the mean
respectively the fluctuating parts of the exact solution and the
approximation provided by the AP-scheme are plotted as a function of
the iteration number. The sequence is initiated with the zero
  function. \modif{With the iterative process, both components converge
  towards the solution until the precision of the schemes is
  reached. At this point, after roughly 27 iterations, the approximation can
  not be improved and a plateau is
  observed.} The convergence of this sequence may be improved thanks to classical relaxation techniques.

\noindent Finally we investigate the positivity of the AP-scheme. With this aim
the anisotropic elliptic problem is solved with a positive source
term, in this case an approximation of the Dirac $\delta$-function. This
function denoted $\delta_a^h$ has a support included in a subset
($[-a,a]\times[-a,a]$, with $0<a<1$) of the simulation domainr
$[-1,1]\times[-1,1]$. Two different parameters $a$ are chosen,
$a=10^{-1}$ and $a=10^{-2}$.

The simulation domain is discretized by a $500\times500$ mesh. For the
smallest value of $a$ the support of the function is reduced to 5
cells in each direction. The source term $\delta_a^h$ is normalized,
such that the maximal value of $\delta_a^h$ grows with vanishing a-parameter.
In table~\ref{table:positiv} the maxima and  minima of the numerical
approximations computed by the AP-scheme ($\phi_A$) and the
discretized Singular Perturbation model ($\phi_{P}$) are gathered for
the two source functions $\delta_a^h$. 
Only large $\eps$-values are considered to verify the positivity of the numerical approximations. Indeed for very small $\varepsilon$ the solution is reduced to its mean part which is the solution of a classical elliptic problem preserving the maximum principle. This means that the relevant question is related to configurations where the fluctuating part $\phi'$ has a significant contribution to the elliptic problem solution. In this range of large and intermediate $\varepsilon$ values, both approximations are comparable. Only slight differences can be observed on the maxima for the smallest $\varepsilon$-parameters. 
The results of table~\ref{table:positiv} demonstrate the positivity of the approximations computed by either the AP-scheme or the  Singular Perturbation model.
\begin{table}[!ht]
\renewcommand\arraystretch{1.25} 
  \centering
\begin{tabular}[c]{|p{0.4cm}|p{1.7cm}|p{1.3cm}|p{1.3cm}|p{1.3cm}|p{1.3cm}|p{1.3cm}|p{1.3cm}|} \cline{2-8} 
 \multicolumn{1}{p{0.4cm}|}{}  &\centering $\varepsilon$ &  $10^2$ & $10$ & $1$ & $10^{-1}$ & $10^{-2}$ & $10^{-3}$\\  \cline{2-8} 
 \end{tabular}

\vspace*{2mm}

\begin{tabular}[c]{|p{0.4cm}|p{1.7cm}|p{1.3cm}|p{1.3cm}|p{1.3cm}|p{1.3cm}|p{1.3cm}|p{1.3cm}|p{1.3cm}|p{1.3cm}|}\hline
\multirow{4}{*}{\rotatebox{90}{$a=10^{-1}$}}
   &max($\phi_{P}$)  & 77.58 & 3.82 & 1.63 & 8.93 & 7.22 & 6.93  \\ \cline{2-8}
   &max($\phi_A$)  & 77.58 & 3.82 & 1.63 & 8.93 & 6.89 & 6.89 \\ \cline{2-8}
   &min($\phi_{P}$)  &$1.9\,10^{-7}$&$2.5\,10^{-7}$&$2.4 \, 10^{-2}$&
   $2.4 \, 10^{-2}$& $2.8 \, 10^{-2}$& $2.8 \, 10^{-2}$ \\ \cline{2-8}
& min($\phi_A$)  &$1.9\,10^{-7}$&$2.5\,10^{-7}$&$2.4 \, 10^{-2}$&
   $2.4 \, 10^{-2}$& $2.8 \, 10^{-2}$& $2.8 \, 10^{-2}$\\ \hline
  \end{tabular}

\vspace*{2mm}

\begin{tabular}[c]{|p{0.4cm}|p{1.7cm}|p{1.3cm}|p{1.3cm}|p{1.3cm}|p{1.3cm}|p{1.3cm}|p{1.3cm}|p{1.3cm}}\hline
\multirow{4}{*}{\rotatebox{90}{$a=10^{-2}$}}
   &max($\phi_{P}$) & $1.8\, 10^{2}$ & $7.1\, 10^{1}$ & $2.6\, 10^{1}$ & $1.2\, 10^{1}$ & $8.29$& $7.34$\\ \cline{2-8}
   &max($\phi_A$) & $1.8\, 10^{2}$ & $7.1\, 10^{1}$ & $2.6\, 10^{1}$ & $1.2\, 10^{1}$ & $7.14$& $7.11$ \\ \cline{2-8}
   &min($\phi_{P}$)  &$1.6 \, 10^{-7}$ &$2.5 \, 10^{-3}$ &$2.4 \, 10^{-2}$ &
   $2.8 \, 10^{-2}$& $2.8 \, 10^{-2}$& $2.8 \, 10^{-2}$\\ \cline{2-8}
& min($\phi_A$) &$1.6 \, 10^{-7}$ &$2.5 \, 10^{-3}$ &$2.4 \, 10^{-2}$ &
   $2.8 \, 10^{-2}$& $2.8 \, 10^{-2}$& $2.8 \, 10^{-2}$\\ \cline{2-8} \hline
  \end{tabular}

  \caption{Maxima and minima of the numerical solutions computed thanks to the AP-scheme ($\phi_A$)  and the Singular Perturbation model ($\phi_{P}$). The elliptic problem is solved with the Dirac $\delta_a^h$ function as a source term on a $500\times500$ mesh.}\label{table:positiv}
\end{table}

\section{Numerical analysis of the AP-scheme}\label{SEC4}
In this last part of the paper we shall concentrate on the numerical
analysis of the $\mathcal{Q}_1$ finite element scheme introduced in
section \ref{SEC31} for solving 
\be \label{E_APi}
\left\{
\begin{array}{l}
\displaystyle - \frac{\partial }{\partial z} \left(A_z \frac{\partial \phi}{\partial z} \right) - \varepsilon \frac{\partial}{\partial x}\left(A_\perp \frac{\partial \phi}{\partial x}  \right) + \varepsilon \frac{\partial}{\partial x}\left(\overline{A_{\perp}' \frac{\partial \phi}{\partial x} }\right) = \varepsilon g\,, \quad \text{in} \quad \Omega\,, \\[3mm]
 \displaystyle \frac{\partial \phi}{\partial z} = 0 \quad \text{on}\quad
\Omega_x \times \partial \Omega_z \,, \qquad
\displaystyle \phi = 0 \quad \text{on} \quad \partial \Omega_x \times \Omega_z\,,
\end{array}
\right.
\ee
where $g \in L^2(\Omega)$ is a given function, with mean value along
the $z$-coordinate equal to zero, $\overline{g}=0$. Moreover we shall
explain why we have to introduce the Lagrange multiplier in order to
solve numerically this equation. We remark that in
contrast to section \ref{SEC3} we omitted for simplicity reasons the
primes for $\phi$, which indicated the fluctuation functions
with zero mean value.\\
The weak form of (\ref{E_APi}) is
\be \label{E_APw}
a(\phi,\psi)=\varepsilon (g,\psi)\,, \quad \forall \psi \in \mathcal{V}\,,
\ee
or equivalently
\be \label{E_APwm}
m(\phi,\psi)=\varepsilon (g,\psi)\,, \quad \forall \psi \in \mathcal{V}\,,
\ee
where $m(\cdot,\cdot)$ is the coercive bilinear form defined in (\ref{E_APmm}).
Let us now consider the corresponding discrete problem
\be \label{E_APad}
a(\phi_h,\psi_h)=\varepsilon (g,\psi_h)\,, \quad \forall \psi_h \in \mathcal{V}_h\,,
\ee
where the finite dimensional space $\mathcal{V}_h \subset \mathcal{V}$ was
introduced in section \ref{SEC31}.
It can be seen that the property $\overline{g}=0$ induces also in the
discrete case that $\overline{\phi_h} = 0$. Thus, following the same arguments as for the continuous case, we can show that equation (\ref{E_APad}) is equivalent to
\be \label{E_APmd}
m(\phi_h,\psi_h)=\varepsilon (g,\psi_h)\,, \quad \forall \psi_h \in \mathcal{V}_h\,.
\ee
The Lax-Milgram theorem implies then the existence and
uniqueness of a discrete solution $\phi_h \in \mathcal{V}_h$. The next theorem gives an estimate of the discretization error $||\phi
-\phi_h||_{\mathcal{V}}$.\\
\modif{ We shall suppose in the sequel, that the diffusion matrices $A_\perp$, $A_z$ and the function $f$ are regular enough, to be able to use standard regularity/interpolation results.}\\
\debthm
Let $\phi \in \mathcal{V}$ be the unique solution of the continuous
problem (\ref{E_APw}) and $\phi_h \in \mathcal{V}_h$ the unique
solution of the discrete
problem (\ref{E_APad}). Both solutions are elements of the \modif{normed} space
$(\mathcal{U},||\cdot||_{\mathcal{U}})$, where
$$
\mathcal{U}:=\{ \psi(\cdot,\cdot) \in \mathcal{V}\,\, / \,\, \overline{\psi}=0
\} \quad \textrm{with}\quad  ||\psi||_{\mathcal{U}}:=||\partial_z \psi||_{L^2(\Omega)}\,.
$$
Then we have the following discretization error estimate
\be \label{E_dis}
||\phi -\phi_h||_{\mathcal{V}}^2=||\partial_z \phi -\partial_z
\phi_h||_{L^2}^2+ \varepsilon ||\partial_x \phi -\partial_x
\phi_h||_{L^2}^2 \le Ch^2\,,
\ee
with a constant $C>0$ independent of $\varepsilon>0$.
Moreover, as $\phi,\phi_h \in \mathcal{U}$, we have 
$$
||\phi -\phi_h||_{\mathcal{U}}^2 \le Ch^2\,.
$$
\finthm
\noindent \debproof
The fact that both solutions $\phi$ and $\phi_h$ belong to the space
$\mathcal{U}$ is an immediate consequence of the fact that the
right-hand  side of the equation (\ref{E_APw}) (resp. (\ref{E_APad}))
satisfies $\overline{g}=0$. The discretization error estimate is
rather standard. Denoting by $\phi_I$ the
interpolant of $\phi$ in the
finite dimensional space $\mathcal{V}_h$, i.e.
$$
\phi_I(x,z):=\sum_{n=1}^{N_x} \sum_{k=1}^{N_z} \phi (x_n,z_k) \chi_n(x)
\kappa_k(z)\,,
$$
we have due to the coercivity of the bilinear form
$m(\cdot,\cdot)$ 
$$
c ||\phi -\phi_h||^2_{\mathcal{V}} \le m(\phi -\phi_h,\phi -\phi_h) =
m(\phi -\phi_h,\phi -\phi_I) \le c||\phi -\phi_h||_{\mathcal{V}}||\phi -\phi_I||_{\mathcal{V}}\,.
$$
Thus
$$
||\phi -\phi_h||_{\mathcal{V}} \le c||\phi -\phi_I||_{\mathcal{V}}\,.
$$
Standard $\mathcal{Q}_1$ finite element interpolation results \modif{\cite{raviart}} yield
for the interpolation error
$$
||\partial_x \phi -\partial_x
\phi_I||_{L^2}^2+  ||\partial_z \phi -\partial_z
\phi_I||_{L^2}^2 \le ch^2 (||\partial_{xx} \phi||_{L^2}^2 +||\partial_{zz} \phi||_{L^2}^2 )\,,
$$
and regularity results for the solution $\phi$ of (\ref{E_APw}), imply $\varepsilon^2  ||\partial_{xx} \phi||_{L^2}^2 +||\partial_{zz}
\phi||_{L^2}^2 \le c \varepsilon^2$. \modif{This last estimate can be found by applying standard $H^2$ regularity results on the solution $\phi_\eps$ of the initial Singular Perturbation problem (\ref{eq:ellipti:original_bis}) (after a change of variable $\xi:=\sqrt{\eps} x$) and then exploiting the decomposition $\phi_\eps=\phi'_\eps+\bar{\phi_\eps}$. Thus, we have altogether} with a constant $c>0$ independent of $\varepsilon >0$
$$
\varepsilon ||\partial_x \phi -\partial_x
\phi_h||_{L^2}^2+  ||\partial_z \phi -\partial_z
\phi_h||_{L^2}^2 \le ch^2\,.
$$
\finproof\\
What is important to observe from the error estimate (\ref{E_dis}) is
that for $\varepsilon \rightarrow 0$ the error $||\phi
-\phi_h||_{H^1}$ in the standard $\varepsilon$-independent $H^1$-norm
blows up. This is one argument why the Singular Perturbation model is inaccurate for
$\varepsilon \ll 1$. However, in the case where $\phi$ and $\phi_h$ are
elements of the space $\mathcal{U}$, we have $||\phi
-\phi_h||_{\mathcal{U}} \le Ch^2$ independently of $\varepsilon$, which means that we have
convergence of the scheme \modif{ in $({\cal U},||\cdot||_{\cal U})$}, uniformly in $\varepsilon>0$. \modif{ The Poincar\'e inequality implies then the uniform convergence in the $||\cdot||_{L^2}$ norm.} The AP-scheme is thus equally
accurate for every value of $0<\varepsilon<1$.\\

\noindent The discretization error $\phi -\phi_h$ is not the only
error we are introducing when solving numerically (\ref{E_APad})
instead of  (\ref{E_APw}). Indeed, (\ref{E_APad}) is nothing but a
linear system
\be \label{E_ls}
M \alpha = v\,,
\ee
to be solved to get the unknowns $\alpha_{nk}:=\phi_h(x_n,z_k)$, where $v_{nk}:=\varepsilon (g, \chi_n
\kappa_k)$ and
the discrete solution of (\ref{E_APad}) is then reconstructed as
$$
\phi_h(x,z)=\sum_{n=1}^{N_x} \sum_{k=1}^{N_z} \alpha_{nk} \chi_n(x)
\kappa_k(z)\,.
$$
Unfortunately the implementation of the system
(\ref{E_ls}) introduces round-off as well as approximation
errors due for example to the numerical computation of $a(\chi_n \kappa_k,\chi_r
\kappa_p)$. Thus the numerical
resolution of (\ref{E_ls}) does not yield the exact solution, but an
approximation $(\tilde{\alpha}_{nk})_{nk}$, solution of the slightly perturbed system
\be \label{E_lsd}
M \tilde{\alpha} = \tilde{v}\,.
\ee
We are now interested in the error estimate $||\phi_h
-\tilde{\phi_h}||_{\mathcal{V}}$, as a function of the perturbation $||v -\tilde{v}||_2$,
where $||\cdot||_2$ denotes the Euclidean norm in $\RR^{N_x N_z}$.\\

\debthm
Let $\alpha$ be the exact solution of (\ref{E_ls}) and
$\tilde{\alpha}$ the exact solution of the perturbed system
(\ref{E_lsd}). Let $\phi_h \in \mathcal{V}_h$ and
$\tilde{\phi_h}\in \mathcal{V}_h$ denote the
corresponding functions
$$
\phi_h(x,z)=\sum_{n=1}^{N_x} \sum_{k=1}^{N_z} \alpha_{nk} \chi_n(x)
\kappa_k(z)\,, \quad \tilde{\phi_h}(x,z)=\sum_{n=1}^{N_x} \sum_{k=1}^{N_z} \tilde{\alpha}_{nk} \chi_n(x)
\kappa_k(z)\,.
$$
Then we have
\be \label{E_round}
 { \varepsilon ||\partial_x \phi_h -\partial_x
\tilde{\phi_h}||_{L^2}^2 +||\partial_z \phi_h -\partial_z
\tilde{\phi_h}||_{L^2}^2 }
\le {c \over \varepsilon}  { ||v-\tilde{v}||_2^2} \,,
\ee
with a constant $c>0$ independent of $\varepsilon>0$ and
$h>0$. However, if both functions $\phi_h$ and $\tilde{\phi_h}$
belong to $\mathcal{U}$, then we have the $\varepsilon$-independent
estimate
$$
||\phi_h -
\tilde{\phi_h}||_{\mathcal{U}} \le c ||v-\tilde{v}||_2 \,.
$$
\finthm

\debproof
Let us denote within this proof $E_{nk}:=
\alpha_{nk}-\tilde{\alpha}_{nk}$ for $n=1,\cdots,N_x$,
$k=1,\cdots,N_z$ and $e_h(x,z):= \phi_h(x,z)
-\tilde{\phi_h}(x,z)$, such that
$$
e_h(x,z)=\sum_{n=1}^{N_x} \sum_{k=1}^{N_z} E_{nk} \chi_n(x)
\kappa_k(z)\,.
$$
Moreover let $N:=N_x N_z$ and $Y \in \RR^{N}$ be an arbitrary vector
associated with the function $y_h(x,z)=\sum_{n=1}^{N_x} \sum_{k=1}^{N_z} Y_{nk} \chi_n(x)
\kappa_k(z)$. Then we have with $(\cdot,\cdot)_2$ the euclidean scalar
product in $\RR^N$ and $M$ the discretization matrix of (\ref{E_ls})
$$
\begin{array}{lll}
||M E||_2 &=&\displaystyle \sup_{Y \in \RR^N\,, Y \neq 0} {(Y,ME)_2 \over ||Y||_2} =
\sup_{Y \in \RR^N\,, Y \neq 0} {m(y_h,e_h) \over ||Y||_2}\,.
\end{array}
$$
Due to the fact that
$$
||Y||_2 \le c ||y_h||_{L^2} \le {c \over \sqrt{\varepsilon}} ||y_h||_{\mathcal{V}}\,,
$$
we have
$$
\begin{array}{lll}
||M E||_2 &\=&\displaystyle \sup_{Y \in \RR^N\,, Y \neq 0} {m(y_h,e_h)
\over ||Y||_2} \ge c \sqrt{\varepsilon} \sup_{y_h \in \mathcal{V}_h \,, y_h\neq 0} {m(y_h,e_h)
\over ||y_h||_{\mathcal{V}}} \ge c \sqrt{\varepsilon} ||e_h||_{\mathcal{V}}\,.
\end{array}
$$
Thus we get with a constant $c>0$ independent of $\varepsilon$
$$
||e_h||_{\mathcal{V}}\le {c \over \sqrt{\varepsilon}} ||M E||_2 = {c \over \sqrt{\varepsilon}} ||v-\tilde{v}||_2\,.
$$
In the case the two functions $\phi_h$ and $\tilde{\phi_h}$
belong to $\mathcal{U}$, i.e. $e_h \in \mathcal{U}$, we can exploit the fact that in
$\mathcal{U}$ \modif{ the Poincar\'e inequality gives rise to}
$||Y||_2 \le c ||y_h||_{L^2} \le c ||y_h||_{\mathcal{U}}$. This yields, as $m(\cdot,\cdot)$ is also coercive on $\mathcal{U}$, that
$$
\begin{array}{lll}
||M E||_2 &\=&\displaystyle \sup_{Y \in \RR^N\,, Y \neq 0} {m(y_h,e_h)
\over ||Y||_2} \ge c  \sup_{y_h \in \mathcal{U} \,, y_h\neq 0} {m(y_h,e_h)
\over ||y_h||_{\mathcal{U}}} \ge c ||e_h||_{\mathcal{U}}\,.
\end{array}
$$
and thus the $\varepsilon$-independent estimate is proved.
\finproof\\
Similarly as for the discretization error, we can deduce from the
round-off error estimate (\ref{E_round}) that, for $\varepsilon \rightarrow 0$, the
standard $H^1$-norm $|| \phi_h -\tilde{\phi_h}||_{H^1}$
explodes. However if we impose that both solutions $\phi_h$ and
$\tilde{\phi_h}$ are elements of the space $\mathcal{U}$, space of
functions with mean value along the $z$-coordinate equal to zero, then
we have the uniform estimate $|| \phi_h -
\tilde{\phi_h}||_{\mathcal{U}} \le c ||v-\tilde{v}||_2 $ \modif{, and by the Poincar\'e inequality $|| \phi_h -
\tilde{\phi_h}||_{L^2} \le c ||v-\tilde{v}||_2$}. Unfortunately even if we know
that $\phi_h \in \mathcal{U}$, this is not necessarily true for
$\tilde{\phi_h}$, if we discretize (\ref{E_APi}). But it can be achieved by forcing the numerical
solution $\tilde{\phi_h}$ to satisfy
$\overline{\tilde{\phi_h}}=0$. Indeed, this can be done by introducing
explicitly in the discrete problem (\ref{E_APad}) the constraint
$\overline{\phi_h}=0$, such that it is much more ingenious to solve instead
\be \label{E_APl}
\left\{
\begin{array}{l}
a(\phi_h,\psi_h)+ b(P_h,\psi_h)=\varepsilon (g,\psi_h)\,,
\quad \forall \psi_h \in \mathcal{V}_h\,, \\[3mm]
b(Q_h,\phi_h)=0\,, \quad \forall Q_h \in \mathcal{W}_h\,,
\end{array}
\right.
\ee
where $\mathcal{W}_h \subset \mathcal{W}$ was constructed in section
\ref{SEC31}. As mentioned in the continuous case this problem is
equivalent for $\eps >0$ to the discrete
problem (\ref{E_APad}). If $\phi_h \in \mathcal{V}_h$ is the unique
solution of (\ref{E_APad}), then $(\phi_h ,0) \in \mathcal{V}_h 
\times \mathcal{W}_h$ is a solution of (\ref{E_APl}). And if $(\phi_h ,P_h) \in \mathcal{V}_h 
\times \mathcal{W}_h$ solves (\ref{E_APl}), then $P_h \equiv
0$ and $\phi_h \in \mathcal{V}_h$ is the unique
solution of (\ref{E_APad}). This last statement is immediately proved
by taking in the variational formulation (\ref{E_APl}) only
$x$-dependent test functions $\psi_h(x)  \in \mathcal{V}_h$. By doing
this, we can be sure that the numerical solution $\tilde{\phi_h}$ of
(\ref{E_APl}) satisfies $\overline{\tilde{\phi_h}}=0$, such that the
error $||\phi_h -\tilde{\phi_h}||_{\mathcal{U}}$ is uniformly
bounded. This proves that the introduction of the constraint
$\overline{\phi_h}=0$ in the AP-formulation is crucial and avoids the
numerical difficulties associated with the original P-model.

\section{Conclusion}\label{SEC5}
In this paper we have introduced an Asymptotic Preserving formulation for
the resolution of a highly anisotropic elliptic equation. We have shown
the advantages of the AP-formulation as compared to the initial
Singular Perturbation  model and to its limit model, when the asymptotic parameter goes to zero. It
came out that the AP-scheme is a powerful tool for the resolution of elliptic problems
presenting huge anisotropies along one coordinate, and gives access to the
simulation in a very easy and precise manner. The
Asymptotic-Preserving method developed here relies on the
decomposition of the solution in its mean part along the anisotropy
direction, and a fluctuation part. This
integration along the anisotropy
direction is easily performed in the context of Cartesian coordinate systems with
one coordinate aligned with the direction of the anisotropy. In a
forthcoming work \cite{brull} this procedure is extended to more
general anisotropies. 

 \section*{Acknowledgments}
This work has been
partially supported by the Marie Curie Actions of the European Commission in the frame of
the DEASE project (MEST-CT-2005-021122) and by the CEA-Cesta in the framework of the 
contracts 'Dynamo-3D' \# 4600108543 and 'Magnefig' \# 06.31.044.

\end{document}